\documentclass[10pt]{article}
\usepackage{latexsym}
\usepackage{amsfonts}
\usepackage{amsmath}
\usepackage{graphicx}
\usepackage[all]{xy}
\usepackage[active]{srcltx} 
\usepackage{hyperref}

\newcommand\tstrut{\rule{0pt}{2.6ex}}

\newcommand{\ba}{\begin{array}}\newcommand{\ea}{\end{array}}

\newcommand{\ns}{\rm}

\newcommand{\nse}{\kern-3pt\ns$=$}\newcommand{\qd}{\hfill$\Box$\medbreak}

\newcommand{\ext}{\raise1pt\hbox{$\ts\bigwedge$}}

\newcommand{\ts}{\textstyle}
\newcommand{\rf}[1]{(\ref{#1})}
\newcommand{\chii}{\raise2pt\hbox{$\chi$}}

\newcommand{\Fg}{\mbox{${\cal F}\kern-2pt_g$}}

\newcommand{\Mg}{\mbox{${\cal M}\kern-2pt_g$}}
\newcommand{\Ng}{\mbox{${\cal N}\kern-2pt_g$}}
\newcommand{\V}{V\kern-1pt}

\newcommand{\Gg}{\mbox{${\cal G}\kern-2pt_g$}}

\renewcommand{\mod}{\mbox{\ns mod\,\,\,}}

\newcommand{\qk}{quaternion-K\"ahler\kern2pt}\renewcommand{\,}{\kern1pt}

\newcommand{\dirac}{/\kern-5pt\partial}

\newcommand{\End}{\rm End}

\newcommand{\lra}{\longrightarrow}
\renewcommand{\ts}{\textstyle}
0

\newtheorem{corol}{Corollary}[section]
\newtheorem{prop}{Proposition}[section]\def\frac#1#2{{#1\over#2}}
\def\be#1\ee{\begin{equation}#1\end{equation}}
\allowdisplaybreaks

\begin{document}
\title{A binary encoding of spinors and applications}

\author{
Gerardo Arizmendi\footnote{Department of Actuarial Sciences, Physics and Mathematics, UDLAP,
San Andr\'es Cholula, Puebla, M\'exico. E-mail: gerardo.arizmendi@udlap.mx} \footnote{Partially supported by a
CONACyT grant
}\,\,\,\,\,\,\,
and
Rafael Herrera\footnote{Centro de Investigaci\'on en Matem\'aticas, A. P. 402,
36000 Guanajuato, Guanajuato, M\'exico. E-mail: rherrera@cimat.mx} $^\dagger$ 
 }

\date{{\ns \today}}

\maketitle

\vspace{-20pt}

{
\abstract{

We present a binary code for spinors and Clifford multiplication using non-negative integers and their binary expressions, which
can be easily implemented in computer programs for explicit calculations.
As applications, we present explicit descriptions of the triality automorphism of $Spin(8)$, explicit representations
of the Lie algebras $\mathfrak{spin}(8)$, $\mathfrak{spin}(7)$ and $\mathfrak{g}_2$, etc. 

}
}


\section{Introduction}

Spinors were first discovered by Cartan in 1913 \cite{Cartan2}, and have been of great relevance in Mathematics and Physics ever since. 
In this note, we introduce a binary code for spinors 
by using a suitable basis and setting up a correspondence between its 
elements and non-negative integers via their binary decompositions.
Such a basis consists of weight vectors of the Spin representation as presented by Friedrich and Sulanke in 1979 
\cite{Friedrich-Sulanke} and which were originally described in a rather different and implicit manner by Brauer and Weyl in 1935
\cite{Brauer-Weyl}. This basis is known to physiscits as a Fock basis \cite{Budinich-Trautman}.
The careful reader will notice that our (integer) encoding uses half the number of bits used by any other binary code of Clifford algebras
\cite{Lounesto-book, Budinich-M}, thus making it computationally more efficient.
Furthermore, Clifford multiplication of vectors with spinors (using the terminology of \cite{friedrich})
becomes a matter of flipping bits in binary expressions and keeping track of powers of $\sqrt{-1}$.
In order to show its usefulness, we develop
very explicit descriptions of some well-known facts such as 
the triality automorphism of $Spin(8)$ (avoiding entirely any reference to the octonions), 
the relationship between Clifford multiplication and the multiplication table of the octonions, and the construction 
of sets of linearly independent vector fields on spheres.
The natural binary/integer code of spinors, as well as the three applications, are the main contributions of the paper. 
The computer implementation of this code has allowed us to carry out many calculations in high dimensions which, otherwise, would
have been impossible. Note that, although we only deal with the 
case of Clifford algebras and spinors defined by positive definite quadratic forms, the binary code can be easily extended to 
semi-definite quadratic forms.

The paper is organized as follows.
In Section \ref{sec: preliminaries}, we recall standard facts about Clifford algebras, the Spin
Lie groups and algebras, the spinor representations and the basis of weight spinors.
In Section \ref{sec: binary}, we introduce the correspondence
between these basic spinors and nonnegative integers, as well as the functions that encode the Clifford multiplication.
In Section \ref{sec: applications}, we present the aforementioned descriptions/applications: triality in dimension 8, the relationship between Clifford multiplication
and the multiplication table of the octonions and, finally, the construction of linearly independent vector fields on spheres.

\section{Preliminaries}\label{sec: preliminaries}

The details of the facts mentioned in this section can be found in \cite{Baum, friedrich}.

\subsection{Clifford algebra, spinors, Clifford multiplication and the Spin group}
Let $Cl_n$ denote the Clifford algebra generated by all the products of the canonical vectors
$e_1, e_2, \ldots, e_n\in\mathbb{R}^n$ subject to the relations
\begin{eqnarray*}
e_j e_k + e_k e_j &=& - \big< e_j, e_k\big> \quad\mbox{for $j\not =k$,} \\
e_j e_j &=& -1 ,
\end{eqnarray*}
where $\big< , \big>$ denotes the standard inner product in $\mathbb{R}^n$, and $Cl_n^0$
the even Clifford subalgebra determined as the fixed point of the involution of $Cl_n$ induced by $-{\rm Id}_{\mathbb{R}^n}$. 
Let
\[\mathbb{C}l_n=Cl_n\otimes_{\mathbb{R}}\mathbb{C},\quad\mbox{and}\quad
\mathbb{C}l_n^0=Cl_n^0\otimes_{\mathbb{R}}\mathbb{C}\]
be the complexifications of $Cl_n$ and $Cl_n^0$. It is well known that
\[\mathbb{C}l_n\cong \left\{
                     \begin{array}{ll}
                     \End(\mathbb{C}^{2^k}), & \mbox{if $n=2k$,}\\
                     \End(\mathbb{C}^{2^k})\oplus\End(\mathbb{C}^{2^k}), & \mbox{if $n=2k+1$,}
                     \end{array},
\right.
\]
where
\[\mathbb{C}^{2^k}=\underbrace{\mathbb{C}^2\otimes \ldots \otimes \mathbb{C}^2}_{\mbox{$k$ times}}\]
is the tensor product of $k=[{n\over 2}]$ copies of $\mathbb{C}^2$.
Let us denote this space as
\[\Delta_n = \mathbb{C}^{2^k},\]
which is called the {\em space of spinors}.
Consider the linear map
\[\kappa:\mathbb{C}l_n \lra \End(\mathbb{C}^{2^k}),\]
which is the aforementioned isomorphism for $n$ even, and the projection onto the first summand for $n$ odd.

The Spin group $Spin(n)\subset Cl_n$ is the subset
\[Spin(n) =\{x_1x_2\cdots x_{2l-1}x_{2l}\,\,|\,\,x_j\in\mathbb{R}^n, \,\,
|x_j|=1,\,\,l\in\mathbb{N}\},\]
endowed with the product of the Clifford algebra.
It is a Lie group and its Lie algebra is
\[\mathfrak{spin}(n)=\mbox{span}\{e_ie_j\,\,|\,\,1\leq i< j \leq n\}.\]
Recall that the Spin group $Spin(n)$ is the universal double cover of $SO(n)$, $n\ge 3$. For $n=2$
we consider $Spin(2)$ to be the connected double cover of $SO(2)$.
The covering map will be denoted by
\[\lambda_n:Spin(n)\rightarrow SO(n),\]
where an element $x_1x_2\cdots x_{2l-1}x_{2l}\in Spin(n)$ is mapped to the orthogonal transformation
\begin{eqnarray*}
\lambda_n: \mathbb{R}^n&\longrightarrow&\mathbb{R}^n\\
 y&\mapsto& x_1x_2\cdots x_{2l-1}x_{2l}\cdot y\cdot x_{2l}x_{2l-1}\cdots x_2x_1.
\end{eqnarray*}
Its differential is given
by
\[\lambda_{n*}(e_ie_j) = 2E_{ij},\]
where $E_{ij}=e_i^*\otimes e_j - e_j^*\otimes e_i$ is the
standard basis of the skew-symmetric matrices and $e^*$ denotes the metric dual of the vector $e$.
We will also denote by $\lambda_n$ the induced representation on
$\ext^*\mathbb{R}^n$.

The restriction of $\kappa$ to $Spin(n)$ defines the Lie group representation
\[
\kappa_n:Spin(n)\lra GL(\Delta_n),\]
which is special unitary. We have the corresponding Lie algebra representation
\[
\kappa_{n*}:\mathfrak{spin}(n)\lra \mathfrak{gl}(\Delta_n),\]
which is, in fact, the restriction of the linear map $\kappa:\mathbb{C}l_n\lra\End(\Delta_n)$ to 
$\mathfrak{spin}(n)\subset \mathbb{C}l_n$.
Note that
$SO(1)=\{1\}$,
$Spin(1)=\{\pm1\}$,
and $\Delta_1 =\mathbb{C}$
the trivial $1$-dimensional representation.

The Clifford multiplication is defined by
\begin{eqnarray*}
\mu_n:\mathbb{R}^n\otimes \Delta_n &\lra&\Delta_n\\
x \otimes \psi &\mapsto& \mu_n(x\otimes \psi)=x\cdot\psi :=\kappa(x)(\psi).
\end{eqnarray*}
It is skew-symmetric with respect to the Hermitian product
\[\left<x\cdot\psi_1 , \psi_2\right>
=-\left<\psi_1 , x\cdot \psi_2\right>,
\]
is $Spin(n)$-equivariant
and  can be extended to a $Spin(n)$-equivariant map
\begin{eqnarray*}
\mu_n:\ext^*(\mathbb{R}^n)\otimes \Delta_n &\lra&\Delta_n\\
\omega \otimes \psi &\mapsto& \omega\cdot\psi.
\end{eqnarray*}
Let
\[{\rm vol}_n  := e_1 e_2\cdots e_n.\]
When $n$ is even, we define the following involution
\begin{eqnarray*}
\Delta_n&\longrightarrow& \Delta_n \\
\psi &\mapsto& (-i)^{n\over 2}{\rm vol}_n\cdot \psi.
\end{eqnarray*}
The $\pm 1$ eigenspaces of this involution are denoted by $\Delta_n^\pm$ and called positive and negative spinors respectively.
These spaces have equal dimension and are irreducible representations of $Spin(n)$
\[\kappa_n^\pm:Spin(n)\lra {\rm Aut}(\Delta_n^\pm).\]
Note that our definition differs from the one given in \cite{friedrich} by a factor $(-1)^{n\over 2}$.

There exist either real or quaternionic structures on the spin representations.
A quaternionic structure $\alpha$ on $\mathbb{C}^2$ is given by
\[\alpha\left(\begin{array}{c}
z_1\\
z_2
              \end{array}
\right) = \left(\begin{array}{c}
-\overline{z}_2\\
\overline{z}_1
              \end{array}\right),\]
and a real structure $\beta$ on $\mathbb{C}^2$ is given by
\[\beta\left(\begin{array}{c}
z_1\\
z_2
              \end{array}
\right) = \left(\begin{array}{c}
\overline{z}_1\\
\overline{z}_2
              \end{array}\right).\]
Note that these structures satisfy
\[
\begin{array}{rclcrcl}
 \left< \alpha(v),w\right> &=& \overline{\left< v,\alpha(w)\right> }, &\quad&
 \left< \alpha(v),\alpha(w)\right> &=& \overline{\left< v,w\right> }, \\
 \left< \beta(v),w\right> &=& \overline{\left< v,\beta(w)\right> }, &\quad&
 \left< \beta(v),\beta(w)\right> &=& \overline{\left< v,w\right> },
\end{array}
\]
with respect to the standard hermitian product in $\mathbb{C}^2$, where $v,w\in \mathbb{C}^2$.
The real and quaternionic structures $\gamma_n$  on $\Delta_n=(\mathbb{C}^2)^{\otimes
[n/2]}$ are built as follows
\[
\begin{array}{cclll}
 \gamma_n &=& (\alpha\otimes\beta)^{\otimes 2k} &\mbox{if $n=8k,8k+1$}& \mbox{(real),} \\
 \gamma_n &=& (\alpha\otimes\beta)^{\otimes 2k}\otimes\alpha &\mbox{if $n=8k+2,8k+3$}&
\mbox{(quaternionic),} \\
 \gamma_n &=& (\alpha\otimes\beta)^{\otimes 2k+1} &\mbox{if $n=8k+4,8k+5$}&\mbox{(quaternionic),} \\
 \gamma_n &=& (\alpha\otimes\beta)^{\otimes 2k+1}\otimes\alpha &\mbox{if $n=8k+6,8k+7$}&\mbox{(real).}
\end{array}
\]
which also satisfy
\[
\begin{array}{rclcrcl}
 \left< \gamma_n(v),w\right> &=& \overline{\left< v,\gamma_n(w)\right> }, &\quad&
 \left< \gamma_n(v),\gamma_n(w)\right> &=& \overline{\left< v,w\right> }, \\
\end{array}
\]
where $v,w\in\Delta_n$.
This means
\[
 \left< v+ \gamma_n(v),w+ \gamma_n(w)\right> \in \mathbb{R}.   \label{eq: real inner product}
\]

Now, we summarize some results about real representations of $Cl_r^0$ in the next table (cf. \cite{Lawson}).
Here $d_r$ denotes the dimension of an irreducible representation of $Cl^0_r$ and $v_r$ the number of distinct
irreducible representations.
\[\begin{array}{|c|c|c|c|c|}
\hline
r \mbox{ (mod 8)}&Cl_r^0&d_r&v_r \tstrut\\
\hline
1&\mathbb R(d_r)&2^{\lfloor\frac r2\rfloor}&1 \tstrut\\
\hline
2&\mathbb C(d_r/2)&2^{\frac r2}&1 \tstrut\\
\hline
3&\mathbb H(d_r/4)&2^{\lfloor\frac r2\rfloor+1}&1 \tstrut\\
\hline
4&\mathbb H(d_r/4)\oplus \mathbb H(d_r/4)&2^{\frac r2}&2 \tstrut\\
\hline
5&\mathbb H(d_r/4)&2^{\lfloor\frac r2\rfloor+1}&1 \tstrut\\
\hline
6&\mathbb C(d_r/2)&2^{\frac r2}&1 \tstrut\\
\hline
7&\mathbb R(d_r)&2^{\lfloor\frac r2\rfloor}&1 \tstrut\\
\hline
8&\mathbb R(d_r)\oplus \mathbb R(d_r)&2^{\frac r2-1}&2 \tstrut\\
\hline
\end{array}
\]
\centerline{Table 1}

Let $\tilde\Delta_r$ denote the real irreducible representation of $Cl_r^0$ for $r\not\equiv0$ $(\mbox{mod } 4) $
and $\tilde\Delta^{\pm}_r$ denote the real irreducible representations for $r\equiv0$ $(\mbox{mod } 4)$. Note that
the representations are complex for $r\equiv 2,6$ $(\mbox{mod } 8)$ and quaternionic for $r\equiv 3,4,5$
$(\mbox{mod } 8)$.

\subsection{A special basis of spinors and an explicit description of $\kappa$ }
In this subsection, we recall the explicit descriptions of $\kappa$ from \cite{friedrich} and the basis of spinors from 
\cite{Friedrich-Sulanke}, which were first discussed by Brauer and Weyl in \cite{Brauer-Weyl}.  

The vectors
\[u_{+1}={1\over \sqrt{2}}(1,-i)\quad\quad\mbox{and}\quad\quad u_{-1}={1\over \sqrt{2}}(1,i),\]
form a unitary basis of $\mathbb{C}^2$.
Consequently, the vectors
\begin{equation}
\{u_{(\varepsilon_1,\ldots,\varepsilon_k)}=u_{\varepsilon_1}\otimes\ldots\otimes
u_{\varepsilon_k}\,\,|\,\, \varepsilon_j=\pm 1,
j=1,\ldots,k\},\label{eq: special basis}
\end{equation}
form a unitary basis of $\Delta_n=(\mathbb{C}^2)^{\otimes
[n/2]}$, which are known to be weight vectors of the Spin representation (see below).

In order to give an explicit description of the map $\kappa$, consider the following matrices with complex entries
\[Id = \left(\begin{array}{ll}
1 & 0\\
0 & 1
       \end{array}\right),\quad
g_1 = \left(\begin{array}{ll}
i & 0\\
0 & -i
       \end{array}\right),\quad
g_2 = \left(\begin{array}{ll}
0 & i\\
i & 0
       \end{array}\right),\quad
T = \left(\begin{array}{ll}
0 & -i\\
i & 0
       \end{array}\right).
\]
Note that
\[g_1(u_{\pm1})= iu_{\mp1},\quad
g_2(u_{\pm1})= \pm u_{\mp1},\quad
T(u_{\pm1})= \mp u_{\pm1}.\]
In general, the generators $e_1, \ldots, e_n$ of the Clifford algebra are mapped under $\kappa$ to the following linear transformations
of $\Delta_n$:
\begin{eqnarray}
 e_1&\mapsto& Id\otimes Id\otimes \ldots\otimes Id\otimes Id\otimes g_1,\nonumber\\
 e_2&\mapsto& Id\otimes Id\otimes \ldots\otimes Id\otimes Id\otimes g_2,\nonumber\\
 e_3&\mapsto& Id\otimes Id\otimes \ldots\otimes Id\otimes g_1\otimes T,\nonumber\\
 e_4&\mapsto& Id\otimes Id\otimes \ldots\otimes Id\otimes g_2\otimes T,\label{eq: explicit Clifford map}\\
 \vdots && \vdots\nonumber\\
 e_{2k-1}&\mapsto& g_1\otimes T\otimes \ldots\otimes T\otimes T\otimes T,\nonumber\\
 e_{2k}&\mapsto& g_2\otimes T\otimes\ldots\otimes T\otimes T\otimes T,\nonumber
\end{eqnarray}
and the last generator
\[ e_{2k+1}\mapsto i\,\, T\otimes T\otimes\ldots\otimes T\otimes T\otimes T\]
if $n=2k+1$. Thus, if $1\leq j\leq k$,
\begin{eqnarray}
 e_{2j-1}u_{\varepsilon_1,\ldots,\varepsilon_k}&=& i(-1)^{j-1}
\left(\prod_{\alpha=k-j+2}^k \varepsilon_{\alpha}\right)
u_{\varepsilon_1,\ldots, (-\varepsilon_{k-j+1}) ,\ldots,\varepsilon_k}    \label{eq: basic clifford multiplication 1}\\
 e_{2j}u_{\varepsilon_1,\ldots,\varepsilon_k}&=&  (-1)^{j-1}
\left(\prod_{\alpha=k-j+1}^k \varepsilon_{\alpha}\right)
u_{\varepsilon_1,\ldots, (-\varepsilon_{k-j+1}) ,\ldots,\varepsilon_k}   \label{eq: basic clifford multiplication 2}
\end{eqnarray}
and
\[
 e_{2k+1}u_{\varepsilon_1,\ldots,\varepsilon_k}= i(-1)^k
\left(\prod_{\alpha=1}^k \varepsilon_{\alpha}\right) u_{\varepsilon_1,\ldots,\varepsilon_k}
\]
if $n=2k+1$ is odd.
Also, the real and quaternionic structures on $\mathbb{C}^2$ look as follows
\begin{eqnarray*}
\alpha(u_{\varepsilon})
   &=& -\varepsilon iu_{-\varepsilon}\\
\beta(u_{\varepsilon})
   &=& u_{-\varepsilon},
\end{eqnarray*}
and, on the basis vectors of $\Delta_n$,
\begin{eqnarray*}
(\alpha\otimes\beta)^{\otimes k}(u_{(\varepsilon_1,...,\varepsilon_{2k})})
   &=&(-i)^k \left(\prod_{s=1}^k \varepsilon_{2s-1}\right) u_{(-\varepsilon_1,...,-\varepsilon_{2k})}\\
\left[(\alpha\otimes\beta)^{\otimes k}\otimes\alpha\right] (u_{(\varepsilon_1,...,\varepsilon_{2k+1})})
   &=&(-i)^{k+1} \left(\prod_{s=1}^{k+1} \varepsilon_{2s-1}\right) u_{(-\varepsilon_1,...,-\varepsilon_{2k+1})}\\
\end{eqnarray*}

{\bf Remark}. From these expressions, we can see that Clifford multiplication by basic vectors amounts to flipping a sign and keeping track of a power of $i=\sqrt{-1}$.

{\bf Remark}.
We became aware of these unitary bases in \cite{Baum} and, since then, we have carried out 
many calculations in many dimensions.
In spite of their usefulness, they are difficult to handle as tensor products in computer programs since they lead to large (sparse) matrices. 
Soon enough, we realized that we could encode such spinors and Clifford multiplication with binary expression (see Section \ref{sec: binary}), and perform calculations in higher dimensions which otherwise would have been impossible.

{\bf Example}. In order to visualize the type of linear transformations given by $\kappa$, 
let us consider $n=6$ and the {\em ordered} spinor basis
\[
\{
u_{(1,1,1)},
u_{(1,1,-1)},
u_{(1,-1,1)},
u_{(1,-1,-1)},
u_{(-1,1,1)},
u_{(-1,1,-1)},
u_{(-1,-1,1)},
u_{(-1,-1,-1)}
\}.
\]
We have, for instance,
\[
\kappa_6(e_1)=\left(\begin{array}{cccccccc}
0 & i &  &  &  &  &  & \\
i & 0 &  &  &  &  &  & \\
 &  & 0 & i &  &  &  & \\
 &  & i & 0 &  &  &  & \\
 &  &  &  & 0 & i &  & \\
 &  &  &  & i & 0 &  & \\
 &  &  &  &  &  & 0 & i\\
 &  &  &  &  &  & i & 0
                  \end{array}
\right),\quad
\kappa_6(e_2)=\left(\begin{array}{cccccccc}
0 & -1 &  &  &  &  &  & \\
1 & 0 &  &  &  &  &  & \\
 &  & 0 & -1 &  &  &  & \\
 &  & 1 & 0 &  &  &  & \\
 &  &  &  & 0 & -1 &  & \\
 &  &  &  & 1 & 0 &  & \\
 &  &  &  &  &  & 0 & -1\\
 &  &  &  &  &  & 1 & 0
                  \end{array}
\right).
\]

\subsubsection{Maximal torus of $Spin(n)$ and weight vectors of $\Delta_n$}

A maximal torus of the group $Spin(n)$ is given by
\[\left\{\prod_{i=1}^{[n/2]} (\cos(\theta_i)+\sin(\theta_i)e_{2i-1}e_{2i})\in Spin(n)\,\,:\,\,
\theta_i\in [0,2\pi], i=1,\ldots,[n/2]\right\}.\]

In order to give a clear idea of the Spin representation, and {\em not for computational purposes}, we will write down some matrices corresponding 
to the transformations of $\Delta_n$ given by $\kappa_n(e_ie_j)$ and $\lambda_n(e_ie_j)$,
for $n=6$ and some $1\leq i<j\leq 6$.
First consider the element $e_1e_2\in Spin(6)$.
In terms of \rf{eq: special basis} and \rf{eq: explicit Clifford map}
\[\kappa_6(e_1e_2)= \left(\begin{array}{cccccccc}
i &  &  &  &  &  &  & \\
 & -i &  &  &  &  &  & \\
 &  & i &  &  &  &  & \\
 &  &  & -i &  &  &  & \\
 &  &  &  & i &  &  & \\
 &  &  &  &  & -i &  & \\
 &  &  &  &  &  & i & \\
 &  &  &  &  &  &  & -i
                      \end{array}
\right),\]
On the other hand, $e_1e_2$ acts on $y\in\mathbb{R}^6$ as follows
\begin{eqnarray*}
e_1e_2(y_1e_1+y_2e_2+y_3e_3+y_4e_4+y_5e_5+y_6e_6)e_2e_1
   &=&
e_1(y_1e_1-y_2e_2+y_3e_3+y_4e_4+y_5e_5+y_6e_6)e_1 \\
   &=&
-y_1e_1-y_2e_2+y_3e_3+y_4e_4+y_5e_5+y_6e_6,
\end{eqnarray*}
i.e.
\[\lambda_6(e_1e_2)=\left(\begin{array}{ccccccc}
-1 &  &  &  &  &   \\
 & -1 &  &  &  &   \\
 &  & 1 &  &  &   \\
 &  &  & 1 &  &   \\
 &  &  &  & 1 &   \\
 &  &  &  &  & 1
        \end{array}
\right)\]
Now, consider the element
\[\cos(\theta)+\sin(\theta)e_1e_2 = (\cos(\theta)e_1+\sin(\theta)e_2)(-e_1)\in Spin(6).\]
On the one hand,
\begin{eqnarray*}
\kappa_6(\cos(\theta)+\sin(\theta)e_1e_2) &=&
\cos(\theta)\left(\begin{array}{cccccccc}
1 &  &  &  &  &  &  & \\
 & 1 &  &  &  &  &  & \\
 &  & 1 &  &  &  &  & \\
 &  &  & 1 &  &  &  & \\
 &  &  &  & 1 &  &  & \\
 &  &  &  &  & 1 &  & \\
 &  &  &  &  &  & 1 & \\
 &  &  &  &  &  &  & 1
                      \end{array}
\right) +
\sin(\theta)\left(\begin{array}{cccccccc}
i &  &  &  &  &  &  & \\
 & -i &  &  &  &  &  & \\
 &  & i &  &  &  &  & \\
 &  &  & -i &  &  &  & \\
 &  &  &  & i &  &  & \\
 &  &  &  &  & -i &  & \\
 &  &  &  &  &  & i & \\
 &  &  &  &  &  &  & -i
                      \end{array}
\right)
\\
&=&
\left(\begin{array}{cccccccc}
e^{i\theta} &  &  &  &  &  &  & \\
 & e^{-i\theta} &  &  &  &  &  & \\
 &  & e^{i\theta} &  &  &  &  & \\
 &  &  & e^{-i\theta} &  &  &  & \\
 &  &  &  & e^{i\theta} &  &  & \\
 &  &  &  &  & e^{-i\theta} &  & \\
 &  &  &  &  &  & e^{i\theta} & \\
 &  &  &  &  &  &  & e^{-i\theta}
                      \end{array}
\right),
\end{eqnarray*}
and, on the other,
\begin{eqnarray*}
(\cos(\theta)e_1+\sin(\theta)e_2)(-e_1)y(-e_1)(\cos(\theta)e_1+\sin(\theta)e_2)
\end{eqnarray*}
induces the transformation on $\mathbb{R}^6$
\[\left(\begin{array}{cccccc}
\cos(2\theta) & -\sin(2\theta) &  &  &  & \\
\sin(2\theta) & \cos(2\theta) &  &  &  & \\
 &  & 1 &  &  & \\
 &  &  & 1 &  & \\
 &  &  &  & 1 & \\
 &  &  &  &  & 1
        \end{array}
\right).\]
Clearly, the two transformations $\kappa_6(\cos(\theta)+\sin(\theta)e_1e_2)$ and $\lambda_6(\cos(\theta)+\sin(\theta)e_1e_2)$ are different.
Setting $\theta=\varphi_1/2$, we see the familiar coefficients $\pm 1/2$ of the Spin representation
\[
\left(\begin{array}{cccccccc}
e^{i\varphi_1\over2} &  &  &  &  &  &  & \\
 & e^{-{i\varphi_1\over2}} &  &  &  &  &  & \\
 &  & e^{i\varphi_1\over2} &  &  &  &  & \\
 &  &  & e^{-{i\varphi_1\over2}} &  &  &  & \\
 &  &  &  & e^{i\varphi_1\over2} &  &  & \\
 &  &  &  &  & e^{-{i\varphi_1\over2}} &  & \\
 &  &  &  &  &  & e^{i\varphi_1\over2} & \\
 &  &  &  &  &  &  & e^{-{i\varphi_1\over2}}
                      \end{array}
\right),
\quad\quad
\left(\begin{array}{cccccc}
\cos(\varphi_1) & -\sin(\varphi_1) &  &  &  & \\
 \sin(\varphi_1) & \cos(\varphi_1) &  &  &  & \\
 &  & 1 &  &  & \\
 &  &  & 1 &  & \\
 &  &  &  & 1 & \\
 &  &  &  &  & 1
        \end{array}
\right).
\]
Similarly,
$\cos(\varphi_2/2)+\sin(\varphi_2/2)e_3e_4$ induces
\[
\left(\begin{array}{cccccccc}
e^{i\varphi_2\over2} &  &  &  &  &  &  & \\
 & e^{{i\varphi_2\over2}} &  &  &  &  &  & \\
 &  & e^{-{i\varphi_2\over2}} &  &  &  &  & \\
 &  &  & e^{-{i\varphi_2\over2}} &  &  &  & \\
 &  &  &  & e^{i\varphi_2\over2} &  &  & \\
 &  &  &  &  & e^{{i\varphi_2\over2}} &  & \\
 &  &  &  &  &  & e^{-{i\varphi_2\over2}} & \\
 &  &  &  &  &  &  & e^{-{i\varphi_2\over2}}
                      \end{array}
\right),
\quad\quad
\left(\begin{array}{cccccc}
1 &  &  &  &  & \\
 & 1 &  &  &  & \\
&&\cos(\varphi_2) & -\sin(\varphi_2) &  &  \\
&&\sin(\varphi_2) & \cos(\varphi_2) &  &   \\
 &  &  &  & 1 & \\
 &  &  &  &  & 1
        \end{array}
\right)
\]
and $\cos(\varphi_3/2)+\sin(\varphi_3/2)e_5e_6$ induces
\[
\left(\begin{array}{cccccccc}
e^{i\varphi_3\over2} &  &  &  &  &  &  & \\
 & e^{{i\varphi_3\over2}} &  &  &  &  &  & \\
 &  & e^{{i\varphi_3\over2}} &  &  &  &  & \\
 &  &  & e^{{i\varphi_3\over2}} &  &  &  & \\
 &  &  &  & e^{-{i\varphi_3\over2}} &  &  & \\
 &  &  &  &  & e^{-{i\varphi_3\over2}} &  & \\
 &  &  &  &  &  & e^{-{i\varphi_3\over2}} & \\
 &  &  &  &  &  &  & e^{-{i\varphi_3\over2}}
                      \end{array}
\right),
\quad\quad
\left(\begin{array}{cccccc}
1 &  &  &  &  & \\
 & 1 &  &  &  & \\
 &  & 1 &  &  & \\
 &  &  & 1 &  & \\
&&&&\cos(\varphi_3) & -\sin(\varphi_3) \\
&&&&\sin(\varphi_3) & \cos(\varphi_3)
        \end{array}
\right)
\]
A general element of the standard maximal torus of $Spin(6)$
\[
{\left(\cos\left({\varphi_1\over2}\right)+{\sin}\left({\varphi_1\over2}\right)e_1e_2\right)}
{\left(\cos\left({\varphi_2\over2}\right)+{\sin}\left({\varphi_2\over2}\right)e_3e_4\right)}
{\left(\cos\left({\varphi_3\over2}\right)+{\sin}\left({\varphi_3\over2}\right)e_5e_6\right)}
\]
has the following matrix representations
\[
\left(\begin{array}{cccccccc}
e^{i(\varphi_1+\varphi_2+\varphi_3)\over2} &  &  &  &  &  &  & \\
 & e^{{i(-\varphi_1+\varphi_2+\varphi_3)\over2}} &  &  &  &  &  & \\
 &  & e^{{i(\varphi_1-\varphi_2+\varphi_3)\over2}} &  &  &  &  & \\
 &  &  & e^{{i(-\varphi_1-\varphi_2+\varphi_3)\over2}} &  &  &  & \\
 &  &  &  & e^{{i(\varphi_1+\varphi_2-\varphi_3)\over2}} &  &  & \\
 &  &  &  &  & e^{{i(-\varphi_1+\varphi_2-\varphi_3)\over2}} &  & \\
 &  &  &  &  &  & e^{{i(\varphi_1-\varphi_2-\varphi_3)\over2}} & \\
 &  &  &  &  &  &  & e^{{i(-\varphi_1-\varphi_2-\varphi_3)\over2}}
                      \end{array}
\right),
\]
\[
\left(\begin{array}{cccccc}
\cos(\varphi_1) & -\sin(\varphi_1) &  &  &  & \\
\sin(\varphi_1) & \cos(\varphi_1) &  &  &  & \\
&&\cos(\varphi_2) & -\sin(\varphi_2) &  &  \\
&&\sin(\varphi_2) & \cos(\varphi_2) &  &   \\
&&&&\cos(\varphi_3) & -\sin(\varphi_3) \\
&&&&\sin(\varphi_3) & \cos(\varphi_3)
        \end{array}
\right).
\]
This formula clearly shows that the basis given in \rf{eq: special basis} is made up of weight vectors
of the spin representation $\Delta_6$.
In general we have
\[\prod_{i=1}^{[n/2]} \left(\cos\left({\varphi_i\over 2}\right)+\sin\left({\varphi_i\over 2}\right)e_{2i-1}e_{2i}\right)\cdot u_{(\varepsilon_1,\ldots,\varepsilon_k)} =
e^{i(\varepsilon_1\varphi_1+\cdots+\varepsilon_k\varphi_k)\over2}u_{(\varepsilon_1,\ldots,\varepsilon_k)}.\]

\section{Binary code}\label{sec: binary}

Given the description in the previous section, we see that the calculation of
$e_{j}u_{\varepsilon_1,\ldots,\varepsilon_k}$, where $k=[n/2]$, depends on $j$, the $k$-tuple
$(\varepsilon_1,\ldots,\varepsilon_k)$ and (possibly on) $n$.
By noticing that
\[+1=(-1)^0 \quad \mbox{and}\quad -1=(-1)^1,\]
we see that for $\varepsilon=\pm1$,
\[\varepsilon=(-1)^{1-\varepsilon\over 2}.\]
Thus, we can change the $k$-tuple $(\varepsilon_1,\ldots,\varepsilon_k)$
by the $k$-tuple
$[{1-\varepsilon_1\over 2},\ldots,{1-\varepsilon_k\over 2}]$ whose entries belong to $\{0,1\}$. Notice that these arrays correspond to the binary expressions of non-negative integers.
For instance, for $n=6$,
\begin{center}
\begin{tabular}[b]{|c|c|c|}\hline
 $(\varepsilon_1,\varepsilon_2,\varepsilon_3)$ & $[{1-\varepsilon_1\over 2},{1-\varepsilon_2\over 2},
 {1-\varepsilon_3\over 2}]$ & Integer \\\hline
(1,1,1) & $[0,0,0] $ & 0\\
(1,1,-1) & $[0,0,1] $ & 1\\
(1,-1,1) & $[0,1,0] $ & 2\\
(1,-1,-1) & $[0,1,1] $& 3\\
(-1,1,1) & $[1,0,0] $& 4\\
(-1,1,-1) & $[1,0,1] $& 5\\
(-1,-1,1) & $[1,1,0] $& 6\\
(-1,-1,-1) & $[1,1,1] $& 7 \\\hline
\end{tabular}
\end{center}
Thus, the aforementioned binary code of spinors is given by the correspondence
\[{(\varepsilon_1,...,\varepsilon_k)}\longleftrightarrow {{1-\varepsilon_k\over 2} (2^0) +{1-\varepsilon_{k-1}\over 2} (2^1) +\ldots
+{1-\varepsilon_2\over 2} (2^{k-2})+{1-\varepsilon_1\over 2} (2^{k-1}) }.\]

{\bf Remark}. The careful reader will notice that for $n=2k$ of $n=2k+1$, this binary encoding of spinors uses $k$ bits 
as opposed to the $2k$ bits of the binary encodings of Clifford algebras \cite{Lounesto-book}.
Since the classical descriptions of the space of spinors are given in terms of minimal ideals within the Clifford algebra itself,  
the inherited binary codes on such minimal ideals use twice as many bits as ours.

\subsection{Clifford multiplication}
The Clifford multiplication of a standard basis vectors $e_p$
with a spinors $u_{a}$, where $a\in\{0,1,2,\ldots,2^{[n/2]}-1\}$ now looks as follows
\begin{eqnarray*}
 e_{2j-1}u_a 
&=& i(-1)^{j-1}(-1)^{\sum_{l=0}^{j-2}[{a/ 2^l}]-2[{a/ 2^{l+1}}]}
 u_{a+(-1)^{[{a/ 2^{j-1}}]-2[{a/ 2^{j}}]}2^{j-1}} \nonumber\\
 e_{2j}u_a 
&=& (-1)^{j-1}(-1)^{\sum_{l=0}^{j-1}[{a/ 2^l}]-2[{a/ 2^{l+1}}]}
 u_{a+(-1)^{[{a/ 2^{j-1}}]-2[{a/ 2^{j}}]}2^{j-1}}, \nonumber
\end{eqnarray*}
 which can be summarized in one formula, for $1\leq p\leq n$, $j:=[{p+1\over 2}]$,
 \begin{eqnarray}
 e_pu_a
 &=& (-1)^{2j-p/2-1+\sum_{l=0}^{j-2}
 a_l+a_{j-1}(-2j+p+1)}u_{a+(-1)^{a_{j-1}}2^{j-1}}\label{formula}
\end{eqnarray}
where
$a_l=\left[{a\over 2^l}\right]-2\left[{a\over 2^{l+1}}\right]$.
Furthermore, if $n$ is odd, $k=[{n\over 2}]$,
\begin{eqnarray}
 e_{2k+1}u_a &=& i(-1)^{k+\sum_{l=0}^{k-1}a_l}
 u_a .
  \label{formula odd}
\end{eqnarray}
These formulas allow us to make general assertions and perform computations in large dimensions without the use of
enormous matrices (recall that the dimension of the $Spin(n)$ representations increases exponentially
with $n$).

{\bf Remark}. This approach also has an important consequence:  
while formulas \rf{eq: basic clifford multiplication 1} and \rf{eq: basic clifford multiplication 2}
seem to depend on $n$, once we write things down using integers
in \rf{formula}, it becomes apparent that Clifford multiplication
does not depend on $n$ if $n$ is even.
For instance, we will always have
\[e_5u_{10}=iu_{15}\]
for all $n\geq 6$. We can actually make the following (non-sharp) claim.

\begin{prop} Let $p\in \mathbb{N}$ and $a\in \mathbb{N}\cup \{0\}$.
Formula \rf{formula} does not depend on $n$ if $n>\max\{p, 1+2\log_2(a+1)\}$.
\end{prop}

In this sense, formula \rf{formula} is rather {\em universal}, but still depends on standard inclusions of Euclidean spaces and their associated Clifford algebras, as well as the explicit mapping of generators \rf{eq: explicit Clifford map}.

\subsubsection{Example: the isomorphism between $\Delta_{2k-1}$ and $\Delta_{2k}^+$}

The space of positive spinors $\Delta_{2k}^+$
is generated by the elements $u_{\varepsilon_1\cdots\varepsilon_k}$ such that
\[\prod_{l=1}^k \varepsilon_l = 1.\]
In the binary code this corresponds to the nonnegative integers whose binary expansion has an even
number of bits $a_l=\left[{a\over 2^l}\right]-2\left[{a\over 2^{l+1}}\right]$ equal to $1$.

Now, the isomorphism 
\[f:\Delta_{2k-1}=\mbox{span}\{u_a\in\mathbb{Z}|0\leq a\leq 2^{k-1}-1 \}
\lra\Delta_{2k}^+=\mbox{span}\left\{u_b\in\mathbb{Z}|0\leq b\leq 2^{k}-1, \sum_{l=0}^{k-1} b_l\equiv 0\,\,(\mod 2) \right\},\]
as representations of the Lie algebra $\mathfrak{spin}(2k-1)$,
is given by
\[f(u_a)= u_{a+\left({1+(-1)^{1+\sum a_l}\over2}\right)2^{k}}.\]
In order to check that the complex linear extension of $f$ is  $\mathfrak{spin}(2k-1)$ equivariant,
let $0\leq a\leq 2^{k-1}$, $1\leq p \leq q\leq 2k-1$. We must verify
\begin{equation}
f(e_pe_q(u_a)) = e_pe_q(f(u_a)). \label{eq: isomorphism}
\end{equation}
Note that the subindices of $u_a$ and $f(u_a)$ have the same binary expression up to the digit corresponding to $2^{k-1}$
so that for $1\leq p< q \leq 2k-2$ the identity \rf{eq: isomorphism}
is fulfilled.
The only cases we have to check are  $1\leq p \leq 2k-2$ and $q=2k-1$.
On the one hand
\[
 e_{2k-1}u_a = i(-1)^{k-1}(-1)^{\sum_{l=0}^{k-2}a_l} u_a
\]
when $u_a$ is considered as a spinor in $\Delta_{2k-1}$
and
\[ e_{2k-1}u_a = i(-1)^{k-1}(-1)^{\sum_{l=0}^{k-2}a_l}
 u_{a+(-1)^{a_{k-1}}2^{k-1}}\]
when $u_a$ is considered as a spinor in $\Delta_{2k}^+$.

On the other hand, in $\Delta_{2k-1}$
\[
e_p e_{2k-1}u_a = i(-1)^{k-1}(-1)^{\sum_{l=0}^{k-2}a_l}
(-1)^{2j-p/2-1}(-1)^{\sum_{l=0}^{j-2}a_l}(-1)^{a_{j-1}(-2j+p+1)}u_{a+(-1)^{a_{j-1}}2^{j-1}}.
\]

{\bf Remark}.
One can even avoid the use of \rf{formula odd} when $n$ is odd and $p=n$ by
using the isomorphism between $\Delta_{2k-1}=\Delta^+_{2k}$.

\section{Applications}\label{sec: applications}

In this long section, we present three applications of the binary code in the form of explicit calculations of the following well-knoun facts: triality in dimension 8 without any reference to the octonions (compare with \cite{Cartan, Harvey,Lounesto,Lounesto-book}), the octonion multiplication table 
(compare with \cite{Baez,Lounesto}) and the construction of independent vector fields on spheres (compare with \cite{Piccinni}).

\subsection{Triality}

We will first recall the idea of triality in a topological form.
As we will recall below, the group $Spin(8)$ is represented
orthogonally on three real 8-dimensional spaces: $\mathbb{R}^8$,
$\tilde\Delta_8^+$ and $\tilde\Delta_8^-$.
In other words, we have three homomorphisms
\[
\xymatrix{
& Spin(8) \ar[d]^{\lambda_8}&\\
& SO(8)&\\
Spin(8)\ar[ur]^{\kappa_8^+}&&Spin(8)\ar[ul]_{\kappa_8^-}\\
}
\]
Now, consider the following two diagrams,
\[
\xymatrix{
  & Spin(8) \ar[d]_{\lambda_8} \\
Spin(8) \ar[ur]^{\sigma}\ar[r]^{\kappa_8^-} & SO(8)
}
\quad\quad\quad
\xymatrix{
  & Spin(8)\ar[d]_{\lambda_8} \\
Spin(8) \ar[ur]^{\tau}\ar[r]^{\kappa_8^+} & SO(8)
}
\]
which include the correspoding lifts $\sigma$ and $\tau$ (due to the simple connectedness of $Spin(8)$).
We will see that $\sigma:Spin(8)\longrightarrow Spin(8)$ is an outer automorphism of order 3 (a {\em triality} automorphism),
$\tau:Spin(8)\longrightarrow Spin(8)$ is an outer automorphism of order 2, and the two automorphisms
generate a copy of the permutation group $S_3$.

First, we will examine the situation explicitly at the Lie algebra level
\begin{equation}
\xymatrix{
  & \mathfrak{spin}(8) \ar[d]^{\lambda_{8*}} \\
\mathfrak{spin}(8) \ar[ur]^{\sigma_*}\ar[r]^{\kappa_{8*}^-} & \mathfrak{so}(8)
}
\quad\quad\quad
\xymatrix{
  & \mathfrak{spin}(8)\ar[d]^{\lambda_{8*}} \\
\mathfrak{spin}(8) \ar[ur]^{\tau_*}\ar[r]^{\kappa_{8*}^+} & \mathfrak{so}(8)
}\label{eq: diagrams Lie algebras}
\end{equation}
and later at the Lie group level.

\subsubsection{The real $Spin(8)$-representations $\tilde\Delta_8^+$ and $\tilde\Delta_8^+$}
Recall that $\gamma_8:\Delta_8\longrightarrow\Delta_8$ is a real structure on $\Delta_8$, which means it is the complexification of a real vector space $\tilde\Delta_8$ given by
\[\tilde\Delta_8 =(1+\gamma_8)(\Delta_8).\]
Furthermore, $\gamma_8$ also preserves the subrepresentations $\Delta_8^+$ and $\Delta_8^-$, i.e. $\gamma_8$ restricts to
real structures on $\Delta_8^+$ and $\Delta_8^-$ and, therefore, they are also complexifications of real vector spaces
\begin{eqnarray*}
\tilde{\Delta}_8^+&=&\{\mbox{$(+1)$-eigenspace of $\gamma_8$ in $\Delta_8^+$} \}\\
&=&(1+\gamma_8)(\Delta_8^+),\\
\tilde{\Delta}_8^-&=&\{\mbox{$(-1)$-eigenspace of $\gamma_8$ in $\Delta_8^-$}\}\\
&=&(1-\gamma_8)(\Delta_8^-).
\end{eqnarray*}
In fact, we have chosen $\tilde\Delta_8^+$ and $\tilde\Delta_8^-$ in this way so that they are compatible with Clifford multiplication
\[\mu_8:\mathbb{R}^8\times\tilde\Delta_8^+\longrightarrow \tilde\Delta_8^-.\]
We have explicit generators for the complex spinor spaces
\begin{eqnarray*}
\Delta_8^+&=&\mathrm{span}(u_0,u_3,u_5,u_6,u_9,u_{10},u_{12},u_{15}),\\
\Delta_8^-&=&\mathrm{span}(u_1,u_2,u_4,u_7,u_8,u_{11},u_{13},u_{14}).
\end{eqnarray*}
For the real representation $\tilde\Delta_8^+$ we have
\begin{eqnarray*}
\tilde\Delta_8^+&=&\mathrm{span}\left\{
   {1\over \sqrt{2}}(u_0-u_{15}),
   {i\over \sqrt{2}}(u_0+u_{15}),
   {1\over \sqrt{2}}(u_3+u_{12}),
   {i\over \sqrt{2}}(u_3-u_{12}),\right.\\
&& \left. {1\over \sqrt{2}}(u_5-u_{10}),
  {i\over \sqrt{2}}(u_5+u_{10}),
  {1\over \sqrt{2}}(u_6+u_9),
  {i\over \sqrt{2}}(u_6- u_9)
\right\}.
\end{eqnarray*}
We choose the ordered basis of $\tilde\Delta_8^-$ to be the image of the basis of $\tilde\Delta_8^+$ under Clifford multiplication by the
canonical vector $e_1\in\mathbb{R} ^8$. Namely,
\begin{eqnarray*}
\tilde\Delta_8^-&=&\mathrm{span}\left\{
 {i\over \sqrt{2}}                      (u_1-u_{14}),
 {-1\over \sqrt{2}}                      (u_1+u_{14}),
  {i\over \sqrt{2}}                     (u_2+u_{13}),
 {1\over \sqrt{2}}                      (u_2- u_{13}),\right.\\
&& \left.{-i\over \sqrt{2}}                     (u_4-u_{11}),
 {-1\over \sqrt{2}}                       (u_4+ u_{11}),
  {i\over \sqrt{2}}                      (u_7+ u_8),
  {-1\over \sqrt{2}}                     (u_7- u_8)
\right\}.
\end{eqnarray*}

\subsubsection{The endomorphism $\sigma_*$}
Using the ordered basis of spinors, one can compute the endomorphisms corresponding
to the generators $e_ie_j\in \mathfrak{spin}(8)$, $1\leq i<j\leq 8$, under the map $\kappa_{8*}^-$ and, in turn, express those endomorphisms as images
of elements of $\mathfrak{spin}(8)$ under $\lambda_{8*}$:
\begin{eqnarray*}
  \kappa_{8*}^-( e_{1} e_{2})&=& -E_{1, 2} - E_{3, 4} - E_{5, 6} - E_{7, 8}
       ={1\over 2}\lambda_{8*}( -e_1e_2 - e_3e_4 - e_5e_6 - e_7e_8),\\
  \kappa_{8*}^-( e_{1} e_{3})&=& -E_{1, 3} + E_{2, 4} - E_{5, 7} + E_{6, 8}
       ={1\over 2}\lambda_{8*}( -e_1e_3 + e_2e_4 - e_5e_7 + e_6e_8),\\
  \kappa_{8*}^-( e_{1} e_{4})&=& -E_{1, 4} - E_{2, 3} + E_{5, 8} + E_{6, 7}
       ={1\over 2}\lambda_{8*}( -e_1e_4 - e_2e_3 + e_5e_8 + e_6e_7),\\
  \kappa_{8*}^-( e_{1} e_{5})&=& -E_{1, 5} + E_{2, 6} + E_{3, 7} - E_{4, 8}
       ={1\over 2}\lambda_{8*}( -e_1e_5 + e_2e_6 + e_3e_7 - e_4e_8),\\
  \kappa_{8*}^-( e_{1} e_{6})&=& -E_{1, 6} - E_{2, 5} - E_{3, 8} - E_{4, 7}
       ={1\over 2}\lambda_{8*}( -e_1e_6 - e_2e_5 - e_3e_8 - e_4e_7),\\
  \kappa_{8*}^-( e_{1} e_{7})&=& -E_{1, 7} + E_{2, 8} - E_{3, 5} + E_{4, 6}
       ={1\over 2}\lambda_{8*}( -e_1e_7 + e_2e_8 - e_3e_5 + e_4e_6),\\
  \kappa_{8*}^-( e_{1} e_{8})&=& -E_{1, 8} - E_{2, 7} + E_{3, 6} + E_{4, 5}
       ={1\over 2}\lambda_{8*}( -e_1e_8 - e_2e_7 + e_3e_6 + e_4e_5),\\
  \kappa_{8*}^-( e_{2} e_{3})&=& E_{1, 4} + E_{2, 3} + E_{5, 8} + E_{6, 7}
       ={1\over 2}\lambda_{8*}( e_1e_4 + e_2e_3 + e_5e_8 + e_6e_7),\\
  \kappa_{8*}^-( e_{2} e_{4})&=& -E_{1, 3} + E_{2, 4} + E_{5, 7} - E_{6, 8}
       ={1\over 2}\lambda_{8*}( e_1e_4 + e_2e_3 + e_5e_8 + e_6e_7),\\
  \kappa_{8*}^-( e_{2} e_{5})&=& E_{1, 6} + E_{2, 5} - E_{3, 8} - E_{4, 7}
       ={1\over 2}\lambda_{8*}( e_1e_6 + e_2e_5 - e_3e_8 - e_4e_7),\\
  \kappa_{8*}^-( e_{2} e_{6})&=& -E_{1, 5} + E_{2, 6} - E_{3, 7} + E_{4, 8}
       ={1\over 2}\lambda_{8*}( -e_1e_5 + e_2e_6 - e_3e_7 + e_4e_8),\\
  \kappa_{8*}^-( e_{2} e_{7})&=& E_{1, 8} + E_{2, 7} + E_{3, 6} + E_{4, 5}
       ={1\over 2}\lambda_{8*}( e_1e_8 + e_2e_7 + e_3e_6 + e_4e_5),\\
  \kappa_{8*}^-( e_{2} e_{8})&=& -E_{1, 7} + E_{2, 8} + E_{3, 5} - E_{4, 6}
       ={1\over 2}\lambda_{8*}( -e_1e_7 + e_2e_8 + e_3e_5 - e_4e_6),\\
  \kappa_{8*}^-( e_{3} e_{4})&=& E_{1, 2} + E_{3, 4} - E_{5, 6} - E_{7, 8}
       ={1\over 2}\lambda_{8*}( e_1e_2 + e_3e_4 - e_5e_6 - e_7e_8),\\
  \kappa_{8*}^-( e_{3} e_{5})&=& E_{1, 7} + E_{2, 8} + E_{3, 5} + E_{4, 6}
       ={1\over 2}\lambda_{8*}( e_1e_7 + e_2e_8 + e_3e_5 + e_4e_6),\\
  \kappa_{8*}^-( e_{3} e_{6})&=& -E_{1, 8} + E_{2, 7} + E_{3, 6} - E_{4, 5}
       ={1\over 2}\lambda_{8*}( -e_1e_8 + e_2e_7 + e_3e_6 - e_4e_5),\\
  \kappa_{8*}^-( e_{3} e_{7})&=& -E_{1, 5} - E_{2, 6} + E_{3, 7} + E_{4, 8}
       ={1\over 2}\lambda_{8*}( -e_1e_5 - e_2e_6 + e_3e_7 + e_4e_8),\\
  \kappa_{8*}^-( e_{3} e_{8})&=& E_{1, 6} - E_{2, 5} + E_{3, 8} - E_{4, 7}
       ={1\over 2}\lambda_{8*}( e_1e_6 - e_2e_5 + e_3e_8 - e_4e_7),\\
  \kappa_{8*}^-( e_{4} e_{5})&=& -E_{1, 8} + E_{2, 7} - E_{3, 6} + E_{4, 5}
       ={1\over 2}\lambda_{8*}( -e_1e_8 + e_2e_7 - e_3e_6 + e_4e_5),\\
  \kappa_{8*}^-( e_{4} e_{6})&=& -E_{1, 7} - E_{2, 8} + E_{3, 5} + E_{4, 6}
       ={1\over 2}\lambda_{8*}( -e_1e_7 - e_2e_8 + e_3e_5 + e_4e_6),\\
  \kappa_{8*}^-( e_{4} e_{7})&=& E_{1, 6} - E_{2, 5} - E_{3, 8} + E_{4, 7}
       ={1\over 2}\lambda_{8*}( e_1e_6 - e_2e_5 - e_3e_8 + e_4e_7),\\
  \kappa_{8*}^-( e_{4} e_{8})&=& E_{1, 5} + E_{2, 6} + E_{3, 7} + E_{4, 8}
       ={1\over 2}\lambda_{8*}( e_1e_5 + e_2e_6 + e_3e_7 + e_4e_8),\\
  \kappa_{8*}^-( e_{5} e_{6})&=& E_{1, 2} - E_{3, 4} + E_{5, 6} - E_{7, 8}
       ={1\over 2}\lambda_{8*}( e_1e_2 - e_3e_4 + e_5e_6 - e_7e_8),\\
  \kappa_{8*}^-( e_{5} e_{7})&=& E_{1, 3} + E_{2, 4} + E_{5, 7} + E_{6, 8}
       ={1\over 2}\lambda_{8*}( e_1e_3 + e_2e_4 + e_5e_7 + e_6e_8),\\
  \kappa_{8*}^-( e_{5} e_{8})&=& -E_{1, 4} + E_{2, 3} + E_{5, 8} - E_{6, 7}
       ={1\over 2}\lambda_{8*}( -e_1e_4 + e_2e_3 + e_5e_8 - e_6e_7),\\
  \kappa_{8*}^-( e_{6} e_{7})&=& -E_{1, 4} + E_{2, 3} - E_{5, 8} + E_{6, 7}
       ={1\over 2}\lambda_{8*}( -e_1e_4 + e_2e_3 - e_5e_8 + e_6e_7),\\
  \kappa_{8*}^-( e_{6} e_{8})&=& -E_{1, 3} - E_{2, 4} + E_{5, 7} + E_{6, 8}
      ={1\over 2}\lambda_{8*}( -e_1e_3 - e_2e_4 + e_5e_7 + e_6e_8),\\
  \kappa_{8*}^-( e_{7} e_{8})&=& E_{1, 2} - E_{3, 4} - E_{5, 6} + E_{7, 8}
      ={1\over 2}\lambda_{8*}( e_1e_2 - e_3e_4 - e_5e_6 + e_7e_8).
\end{eqnarray*}
This means, in terms of the first diagram in \rf{eq: diagrams Lie algebras},
\begin{eqnarray*}
 \sigma_*( e_{1} e_{2})&=&{1\over 2} ( -e_1e_2 - e_3e_4 - e_5e_6 - e_7e_8),\\
 \sigma_*( e_{1} e_{3})&=&{1\over 2} ( -e_1e_3 + e_2e_4 - e_5e_7 + e_6e_8),\\
 \sigma_*( e_{1} e_{4})&=&{1\over 2} ( -e_1e_4 - e_2e_3 + e_5e_8 + e_6e_7),\\
 \sigma_*( e_{1} e_{5})&=&{1\over 2} ( -e_1e_5 + e_2e_6 + e_3e_7 - e_4e_8),\\
 \sigma_*( e_{1} e_{6})&=&{1\over 2} ( -e_1e_6 - e_2e_5 - e_3e_8 - e_4e_7),\\
 \sigma_*( e_{1} e_{7})&=&{1\over 2} ( -e_1e_7 + e_2e_8 - e_3e_5 + e_4e_6),\\
 \sigma_*( e_{1} e_{8})&=&{1\over 2} ( -e_1e_8 - e_2e_7 + e_3e_6 + e_4e_5),\\
 \sigma_*( e_{2} e_{3})&=&{1\over 2} ( e_1e_4 + e_2e_3 + e_5e_8 + e_6e_7),\\
 \sigma_*( e_{2} e_{4})&=&{1\over 2} ( -e_1e_3 + e_2e_4 + e_5e_7 - e_6e_8),\\
 \sigma_*( e_{2} e_{5})&=&{1\over 2} ( e_1e_6 + e_2e_5 - e_3e_8 - e_4e_7),\\
 \sigma_*( e_{2} e_{6})&=&{1\over 2} ( -e_1e_5 + e_2e_6 - e_3e_7 + e_4e_8),\\
 \sigma_*( e_{2} e_{7})&=&{1\over 2} ( e_1e_8 + e_2e_7 + e_3e_6 + e_4e_5),\\
 \sigma_*( e_{2} e_{8})&=&{1\over 2} ( -e_1e_7 + e_2e_8 + e_3e_5 - e_4e_6),\\
 \sigma_*( e_{3} e_{4})&=&{1\over 2} ( e_1e_2 + e_3e_4 - e_5e_6 - e_7e_8),\\
 \sigma_*( e_{3} e_{5})&=&{1\over 2} ( e_1e_7 + e_2e_8 + e_3e_5 + e_4e_6),\\
 \sigma_*( e_{3} e_{6})&=&{1\over 2} ( -e_1e_8 + e_2e_7 + e_3e_6 - e_4e_5),\\
 \sigma_*( e_{3} e_{7})&=&{1\over 2} ( -e_1e_5 - e_2e_6 + e_3e_7 + e_4e_8),\\
 \sigma_*( e_{3} e_{8})&=&{1\over 2} ( e_1e_6 - e_2e_5 + e_3e_8 - e_4e_7),\\
 \sigma_*( e_{4} e_{5})&=&{1\over 2} ( -e_1e_8 + e_2e_7 - e_3e_6 + e_4e_5),\\
 \sigma_*( e_{4} e_{6})&=&{1\over 2} ( -e_1e_7 - e_2e_8 + e_3e_5 + e_4e_6),\\
 \sigma_*( e_{4} e_{7})&=&{1\over 2} ( e_1e_6 - e_2e_5 - e_3e_8 + e_4e_7),\\
 \sigma_*( e_{4} e_{8})&=&{1\over 2} ( e_1e_5 + e_2e_6 + e_3e_7 + e_4e_8),\\
 \sigma_*( e_{5} e_{6})&=&{1\over 2} ( e_1e_2 - e_3e_4 + e_5e_6 - e_7e_8),\\
 \sigma_*( e_{5} e_{7})&=&{1\over 2} ( e_1e_3 + e_2e_4 + e_5e_7 + e_6e_8),\\
 \sigma_*( e_{5} e_{8})&=&{1\over 2} ( -e_1e_4 + e_2e_3 + e_5e_8 - e_6e_7),\\
 \sigma_*( e_{6} e_{7})&=&{1\over 2} ( -e_1e_4 + e_2e_3 - e_5e_8 + e_6e_7),\\
 \sigma_*( e_{6} e_{8})&=&{1\over 2} ( -e_1e_3 - e_2e_4 + e_5e_7 + e_6e_8),\\
 \sigma_*( e_{7} e_{8})&=&{1\over 2} ( e_1e_2 - e_3e_4 - e_5e_6 + e_7e_8).
\end{eqnarray*}
In other words, we have defined $\sigma_*$ in such a way that
\[\lambda_{8*}\circ \sigma_* = \kappa_{8*}^-.\]
In order to show that $\sigma_*$ is of order 3, let us consider, for instance
\[ \sigma_*( e_{1} e_{2})={1\over 2} ( -e_1e_2 - e_3e_4 - e_5e_6 - e_7e_8).
\]
Then
\begin{eqnarray*}
  \sigma_*(\sigma_*( e_{1} e_{2}))
     &=&{1\over 2} ( -\sigma_*(e_1e_2) - \sigma_*(e_3e_4) - \sigma_*(e_5e_6) - \sigma_*(e_7e_8))\\
     &=&{1\over 4} ( -( -e_1e_2 - e_3e_4 - e_5e_6 - e_7e_8) - ( e_1e_2 + e_3e_4 - e_5e_6 - e_7e_8) \\
     & &   - ( e_1e_2 - e_3e_4 + e_5e_6 - e_7e_8) - ( e_1e_2 - e_3e_4 - e_5e_6 + e_7e_8))\\
     &=&{1\over 2} (-e_1e_2 +e_3e_4 +e_5e_6+e_7e_8),
\end{eqnarray*}
and
\begin{eqnarray*}
  \sigma_*(\sigma_*(\sigma_*( e_{1} e_{2})))
     &=&{1\over 2} (-\sigma_*(e_1e_2) +\sigma_*(e_3e_4) +\sigma_*(e_5e_6)+\sigma_*(e_7e_8))\\
     &=&{1\over 4} (-( -e_1e_2 - e_3e_4 - e_5e_6 - e_7e_8) +( e_1e_2 + e_3e_4 - e_5e_6 - e_7e_8) \\
     & & +( e_1e_2 - e_3e_4 + e_5e_6 - e_7e_8)+( e_1e_2 - e_3e_4 - e_5e_6 + e_7e_8))\\
     &=&e_1e_2 .
\end{eqnarray*}
All the other cases are similar.
In fact, using the standard ordered basis $\{e_1e_2,e_1e_3,...,e_7e_8\}$
of $\mathfrak{spin}(8)$ we have the matrix representation
{\tiny
\[\sigma_*={1\over 2}
\left(
\begin{array}{cccccccccccccccccccccccccccc}
     -1  & 0 & 0 & 0 & 0 & 0 & 0 & 0 & 0 & 0 & 0 & 0 & 0 & 1  & 0 &
    0 & 0 & 0 & 0 & 0 & 0 & 0 & 1  & 0 & 0 & 0 & 0 & 1 \\
     0 & -1  & 0 & 0 & 0 & 0 & 0 & 0 & -1  & 0 & 0 & 0 & 0 & 0 & 0
    & 0 & 0 & 0 & 0 & 0 & 0 & 0 & 0 & 1  & 0 & 0 & -1  & 0\\
     0 & 0 & -1  & 0 & 0 & 0 & 0 & 1  & 0 & 0 & 0 & 0 & 0 & 0 & 0 &
    0 & 0 & 0 & 0 & 0 & 0 & 0 & 0 & 0 & -1  & -1  & 0 & 0\\
     0 & 0 & 0 & -1  & 0 & 0 & 0 & 0 & 0 & 0 & -1  & 0 & 0 & 0 & 0
    & 0 & -1  & 0 & 0 & 0 & 0 & 1  & 0 & 0 & 0 & 0 & 0 & 0\\
     0 & 0 & 0 & 0 & -1  & 0 & 0 & 0 & 0 & 1  & 0 & 0 & 0 & 0 & 0 &
    0 & 0 & 1  & 0 & 0 & 1  & 0 & 0 & 0 & 0 & 0 & 0 & 0\\
     0 & 0 & 0 & 0 & 0 & -1  & 0 & 0 & 0 & 0 & 0 & 0 & -1  & 0 &
    1  & 0 & 0 & 0 & 0 & -1  & 0 & 0 & 0 & 0 & 0 & 0 & 0 & 0\\
     0 & 0 & 0 & 0 & 0 & 0 & -1  & 0 & 0 & 0 & 0 & 1  & 0 & 0 & 0 &
    -1  & 0 & 0 & -1  & 0 & 0 & 0 & 0 & 0 & 0 & 0 & 0 & 0\\
     0 & 0 & -1  & 0 & 0 & 0 & 0 & 1  & 0 & 0 & 0 & 0 & 0 & 0 & 0 &
    0 & 0 & 0 & 0 & 0 & 0 & 0 & 0 & 0 & 1  & 1  & 0 & 0\\
     0 & 1  & 0 & 0 & 0 & 0 & 0 & 0 & 1  & 0 & 0 & 0 & 0 & 0 & 0 &
    0 & 0 & 0 & 0 & 0 & 0 & 0 & 0 & 1  & 0 & 0 & -1  & 0\\
     0 & 0 & 0 & 0 & -1  & 0 & 0 & 0 & 0 & 1  & 0 & 0 & 0 & 0 & 0 &
    0 & 0 & -1  & 0 & 0 & -1  & 0 & 0 & 0 & 0 & 0 & 0 & 0\\
     0 & 0 & 0 & 1  & 0 & 0 & 0 & 0 & 0 & 0 & 1  & 0 & 0 & 0 & 0 &
    0 & -1  & 0 & 0 & 0 & 0 & 1  & 0 & 0 & 0 & 0 & 0 & 0\\
     0 & 0 & 0 & 0 & 0 & 0 & -1  & 0 & 0 & 0 & 0 & 1  & 0 & 0 & 0 &
    1  & 0 & 0 & 1  & 0 & 0 & 0 & 0 & 0 & 0 & 0 & 0 & 0\\
     0 & 0 & 0 & 0 & 0 & 1  & 0 & 0 & 0 & 0 & 0 & 0 & 1  & 0 & 1
    & 0 & 0 & 0 & 0 & -1  & 0 & 0 & 0 & 0 & 0 & 0 & 0 & 0\\
     -1  & 0 & 0 & 0 & 0 & 0 & 0 & 0 & 0 & 0 & 0 & 0 & 0 & 1  & 0 &
    0 & 0 & 0 & 0 & 0 & 0 & 0 & -1  & 0 & 0 & 0 & 0 & -1 \\
     0 & 0 & 0 & 0 & 0 & -1  & 0 & 0 & 0 & 0 & 0 & 0 & 1  & 0 & 1
    & 0 & 0 & 0 & 0 & 1  & 0 & 0 & 0 & 0 & 0 & 0 & 0 & 0\\
     0 & 0 & 0 & 0 & 0 & 0 & 1  & 0 & 0 & 0 & 0 & 1  & 0 & 0 & 0 &
    1  & 0 & 0 & -1  & 0 & 0 & 0 & 0 & 0 & 0 & 0 & 0 & 0\\
     0 & 0 & 0 & 1  & 0 & 0 & 0 & 0 & 0 & 0 & -1  & 0 & 0 & 0 & 0 &
    0 & 1  & 0 & 0 & 0 & 0 & 1  & 0 & 0 & 0 & 0 & 0 & 0\\
     0 & 0 & 0 & 0 & -1  & 0 & 0 & 0 & 0 & -1  & 0 & 0 & 0 & 0 & 0
    & 0 & 0 & 1  & 0 & 0 & -1  & 0 & 0 & 0 & 0 & 0 & 0 & 0\\
     0 & 0 & 0 & 0 & 0 & 0 & 1  & 0 & 0 & 0 & 0 & 1  & 0 & 0 & 0 &
    -1  & 0 & 0 & 1  & 0 & 0 & 0 & 0 & 0 & 0 & 0 & 0 & 0\\
     0 & 0 & 0 & 0 & 0 & 1  & 0 & 0 & 0 & 0 & 0 & 0 & -1  & 0 & 1
    & 0 & 0 & 0 & 0 & 1  & 0 & 0 & 0 & 0 & 0 & 0 & 0 & 0\\
     0 & 0 & 0 & 0 & -1  & 0 & 0 & 0 & 0 & -1  & 0 & 0 & 0 & 0 & 0
    & 0 & 0 & -1  & 0 & 0 & 1  & 0 & 0 & 0 & 0 & 0 & 0 & 0\\
     0 & 0 & 0 & -1  & 0 & 0 & 0 & 0 & 0 & 0 & 1  & 0 & 0 & 0 & 0 &
    0 & 1  & 0 & 0 & 0 & 0 & 1  & 0 & 0 & 0 & 0 & 0 & 0\\
     -1  & 0 & 0 & 0 & 0 & 0 & 0 & 0 & 0 & 0 & 0 & 0 & 0 & -1  & 0
    & 0 & 0 & 0 & 0 & 0 & 0 & 0 & 1  & 0 & 0 & 0 & 0 & -1 \\
     0 & -1  & 0 & 0 & 0 & 0 & 0 & 0 & 1  & 0 & 0 & 0 & 0 & 0 & 0 &
    0 & 0 & 0 & 0 & 0 & 0 & 0 & 0 & 1  & 0 & 0 & 1  & 0\\
     0 & 0 & 1  & 0 & 0 & 0 & 0 & 1  & 0 & 0 & 0 & 0 & 0 & 0 & 0 &
    0 & 0 & 0 & 0 & 0 & 0 & 0 & 0 & 0 & 1  & -1  & 0 & 0\\
     0 & 0 & 1  & 0 & 0 & 0 & 0 & 1  & 0 & 0 & 0 & 0 & 0 & 0 & 0 &
    0 & 0 & 0 & 0 & 0 & 0 & 0 & 0 & 0 & -1  & 1  & 0 & 0\\
     0 & 1  & 0 & 0 & 0 & 0 & 0 & 0 & -1  & 0 & 0 & 0 & 0 & 0 & 0 &
    0 & 0 & 0 & 0 & 0 & 0 & 0 & 0 & 1  & 0 & 0 & 1  & 0\\
     -1  & 0 & 0 & 0 & 0 & 0 & 0 & 0 & 0 & 0 & 0 & 0 & 0 & -1  & 0
    & 0 & 0 & 0 & 0 & 0 & 0 & 0 & -1  & 0 & 0 & 0 & 0 & 1 \\
\end{array}
\right),
\]
}
which one can verify is of order 3.
Furthermore, the map $\sigma_*$ has eigenvalues
\[e^{2\pi i\over3},e^{-2\pi i\over3},1,\]
with multiplicities 7, 7 and 14 respectively.
The eigenspace corresponding to $1$ is generated by
\begin{eqnarray}
&&                    \left\{ 
                     e_{2} e_{3} + e_{6} e_{7},
                     -e_{2} e_{4} + e_{6} e_{8},
                     -e_{3} e_{4} + e_{7} e_{8},
                     -e_{2} e_{6} + e_{3} e_{7},
                     -e_{2} e_{5} + e_{3} e_{8},                     
                     e_{2} e_{7} + e_{4} e_{5},
                     -e_{2} e_{8} + e_{4} e_{6},
                     \right.\nonumber\\
&&\left.                     
                     -e_{2} e_{5} + e_{4} e_{7},
                     e_{2} e_{6} + e_{4} e_{8},
                     -e_{3} e_{4} + e_{5} e_{6},
                     e_{2} e_{4} + e_{5} e_{7},
                     e_{2} e_{3} + e_{5} e_{8},
                     e_{2} e_{8} + e_{3} e_{5},
                     e_{2} e_{7} + e_{3} e_{6},
                     \right\}\label{eq: generators g2}
\end{eqnarray}
which generates a copy of $\mathfrak{g}_2\subset \mathfrak{spin}(8)\subset Cl_8^0$ 
(see \cite{Baez,Harvey,Lounesto} for definitions of $\mathfrak{g}_2$ and $G_2$).
Note that none of the generators includes the vector $e_1$, which makes this copy of $\mathfrak{g}_2$
a subalgebra of the copy of $\mathfrak{spin}(7)$ generated by the span of $\{e_2,e_3,e_4,e_5,e_6,e_7,e_8\}$.
This copy of $\mathfrak{g}_2$ annihilates the basic positive spinor
\[{1\over \sqrt{2}}(u_0-u_{15})\in \tilde\Delta_8^+,\]
so that $\tilde\Delta_8^+=\mathbf{1}\oplus\mathbb{R}^7$ under $\mathfrak{g}_2$ and also
annihilates the basic negative spinor
\[{i\over \sqrt{2}}                      (u_1-u_{14})\in \tilde\Delta_8^-\]
under Clifford multiplication, so that $\tilde\Delta_8^-=\mathbf{1}\oplus\mathbb{R}^7$ under $\mathfrak{g}_2$.
The matrix representation
for a general element
\begin{eqnarray*}
&&                    
                     \alpha_{1} (e_{2} e_{3} + e_{6} e_{7})+
                     \alpha_{2} (-e_{2} e_{4} + e_{6} e_{8})+
                     \alpha_{3} (-e_{3} e_{4} + e_{7} e_{8})+
                     \alpha_{4} (-e_{2} e_{6} + e_{3} e_{7})+
                     \alpha_{5} (-e_{2} e_{5} + e_{3} e_{8})\\
                     &&+                     
                     \alpha_{6} (e_{2} e_{7} + e_{4} e_{5})+
                     \alpha_{7} (-e_{2} e_{8} + e_{4} e_{6})+
                     \alpha_{8} (-e_{2} e_{5} + e_{4} e_{7})+
                     \alpha_{9} (e_{2} e_{6} + e_{4} e_{8})+
                     \alpha_{10} (-e_{3} e_{4} + e_{5} e_{6})\\
                     &&+
                     \alpha_{11} (e_{2} e_{4} + e_{5} e_{7})+
                     \alpha_{12} (e_{2} e_{3} + e_{5} e_{8})+
                     \alpha_{13} (e_{2} e_{8} + e_{3} e_{5})+
                     \alpha_{14} (e_{2} e_{7} + e_{3} e_{6})                    
\end{eqnarray*}
on both $\tilde\Delta_8^+$ and $\tilde\Delta_8^-$, is
\[
2\left(
\begin{array}{cccccccc}
0&0&0&0&0&0&0&0\\
0&0&  \alpha_{1}+  \alpha_{12}&  -\alpha_{2}+ \alpha_{11}&-\alpha_{5}-  \alpha_{8}&  -\alpha_{4}+  \alpha_{9}&  \alpha_{6}+  \alpha_{14}&-  \alpha_{7}+  \alpha_{13}\\
0&-  \alpha_{1}-  \alpha_{12}&0&-  \alpha_{3}-  \alpha_{10}&  \alpha_{13}&  \alpha_{14}&  \alpha_{4 }&  \alpha_{5}\\
0&  \alpha_{2}-  \alpha_{11}&  \alpha_{3}+  \alpha_{10}&0&  \alpha_{6}&  \alpha_{7}&  \alpha_{8 }&  \alpha_{9}\\
0&  \alpha_{5 }+  \alpha_{8}&-  \alpha_{13}&-  \alpha_{6}&0&  \alpha_{10}&  \alpha_{11}&  \alpha_{12}\\
0&  \alpha_{4}-  \alpha_{9 }&-  \alpha_{14}&-  \alpha_{7}&-  \alpha_{10}&0&  \alpha_{1}&  \alpha_{2}\\
0&-  \alpha_{6}-  \alpha_{14}&-  \alpha_{4 }&-  \alpha_{8 }&-  \alpha_{11}&- \alpha_{1}&0& \alpha_{3}\\
0&  \alpha_{7}-  \alpha_{13}&-  \alpha_{5}&-  \alpha_{9}&-  \alpha_{12}&-  \alpha_{2}&- \alpha_{3}&0\\
\end{array}
\right).
\]
Let us compute an explicit element of the group $G_2$.
Consider the element $e_{2} e_{3} + e_{6} e_{7}\in\mathfrak{g}_2\subset \mathfrak{spin}(8)\subset Cl_8^0$, 
and the one parameter subgroup
\begin{eqnarray*}
exp(t(e_{2} e_{3} + e_{6} e_{7}))
   &=& exp(te_{2} e_{3}) exp(te_{6} e_{7})\\
   &=& (\cos(t)+\sin(t)e_2e_3)(\cos(t)+\sin(t)e_6e_7)\\
   &=& {1\over 2}(-e_2e_3e_6e_7\cos(2t)+e_2e_3e_6e_7+e_2e_3\sin(2t)+e_6e_7\sin(2t)+\cos(2t)+1)\\
   &=& {1\over 2}(e_2e_3e_6e_7(1-\cos(2t))+(e_2e_3+e_6e_7)\sin(2t)+\cos(2t)+1)\\
&\in&    G_2\subset Spin(7)\subset Spin(8)\subset Cl_8^0.
\end{eqnarray*} 
Its image under $\kappa_8^-$ is
\begin{eqnarray*}
\kappa_8^-(exp(t(e_{2} e_{3} + e_{6} e_{7})))
   &=& (\cos(t){\rm Id}_8+\sin(t)\kappa_8^-(e_2e_3))(\cos(t){\rm Id}_8+\sin(t)\kappa_8^-(e_6e_7))\\
   &=&\left(\begin{array}{cccccccc}
   1 & • & • & • & • & • & • & • \\ 
   • & \cos(2t) & -\sin(2t) & • & • & • & • & • \\ 
   • & \sin(2t) & \cos(2t) & • & • & • & • & • \\ 
   • & • & • & 1 & • & • & • & • \\ 
   • & • & • & • & 1 & • & • & • \\ 
   • & • & • & • & • & \cos(2t) & -\sin(2t) & • \\ 
   • & • & • & • & • & \sin(2t) & \cos(2t) & • \\ 
   • & • & • & • & • & • & • & 1
   \end{array} \right)\\
&\in &    \kappa_8^-( G_2)\subset SO(7)\subset SO(8).
\end{eqnarray*}

The eigenspace corresponding to $e^{2\pi i\over3}$ is generated by
\begin{eqnarray*}
    && \{e_{6} e_{8} - e_{5} e_{7} + e_{2} e_{4} + i e_{1} e_{3} \sqrt{3},
     e_{4} e_{7} + e_{3} e_{8} + e_{2} e_{5} - i e_{1} e_{6} \sqrt{3}     ,
     e_{4} e_{6} - e_{3} e_{5} + e_{2} e_{8} + i e_{1} e_{7} \sqrt{3}     ,\\
&&     e_{7} e_{8} + e_{5} e_{6} + e_{3} e_{4} - i e_{1} e_{2} \sqrt{3}     ,
     e_{4} e_{5} + e_{3} e_{6} - e_{2} e_{7} + i e_{1} e_{8} \sqrt{3}     ,
     e_{4} e_{8} - e_{3} e_{7} - e_{2} e_{6} - i e_{1} e_{5} \sqrt{3}     ,\\
&&     e_{6} e_{7} + e_{5} e_{8} - e_{2} e_{3} + i e_{1} e_{4} \sqrt{3}    \},
\end{eqnarray*}
and eigenspace corresponding to $e^{-2\pi i\over3}$ is generated by
\begin{eqnarray*}
&&\{     e_{6} e_{8} - e_{5} e_{7} + e_{2} e_{4} - i e_{1} e_{3} \sqrt{3}     ,
     e_{4} e_{7} + e_{3} e_{8} + e_{2} e_{5} + i e_{1} e_{6} \sqrt{3}    ,
     e_{4} e_{6} - e_{3} e_{5} + e_{2} e_{8} - i e_{1} e_{7} \sqrt{3}     ,\\
 &&    e_{7} e_{8} + e_{5} e_{6} + e_{3} e_{4} + i e_{1} e_{2} \sqrt{3}    , 
     e_{4} e_{5} + e_{3} e_{6} - e_{2} e_{7} - i e_{1} e_{8} \sqrt{3}     ,
     e_{4} e_{8} - e_{3} e_{7} - e_{2} e_{6} + i e_{1} e_{5} \sqrt{3}     ,\\
&&     e_{6} e_{7} + e_{5} e_{8} - e_{2} e_{3} - i e_{1} e_{4} \sqrt{3}     \}
\end{eqnarray*}

\subsubsection{The endomorphism $\tau_*$}
Using the ordered basis of spinors, one can compute the endomorphisms corresponding
to the elements $e_ie_j\in \mathfrak{spin}(8)$, $1\leq i< j \leq 8$, under the map $\kappa_{8*}^+$ and, in turn, express those endomorphisms as images
of elements of $\mathfrak{spin}(8)$ under $\lambda_{8*}$:
\begin{eqnarray*}
 \kappa_{8*}^+( e_1 e_2)
   &=& E_{1, 2} + E_{3, 4} + E_{5, 6} + E_{7, 8}
    ={1\over 2}\lambda_{8*}(e_{1, 2} + e_{3, 4} + e_{5, 6} + e_{7, 8}),\\
 \kappa_{8*}^+( e_1 e_3)
    &=& E_{1, 3} - E_{2, 4} + E_{5, 7} - E_{6, 8}
    ={1\over 2}\lambda_{8*}( e_{1, 3} - e_{2, 4} + e_{5, 7} - e_{6, 8}),\\
 \kappa_{8*}^+( e_1 e_4)
    &=& E_{1, 4} + E_{2, 3} - E_{5, 8} - E_{6, 7}
    ={1\over 2}\lambda_{8*}( e_{1, 4} + e_{2, 3} - e_{5, 8} - e_{6, 7}),\\
 \kappa_{8*}^+( e_1 e_5)
    &=& E_{1, 5} - E_{2, 6} - E_{3, 7} + E_{4, 8}
    ={1\over 2}\lambda_{8*}( e_{1, 5} - e_{2, 6} - e_{3, 7} + e_{4, 8}),\\
 \kappa_{8*}^+( e_1 e_6)
    &=& E_{1, 6} + E_{2, 5} + E_{3, 8} + E_{4, 7}
    ={1\over 2}\lambda_{8*}( e_{1, 6} + e_{2, 5} + e_{3, 8} + e_{4, 7}),\\
 \kappa_{8*}^+( e_1 e_7)
    &=& E_{1, 7} - E_{2, 8} + E_{3, 5} - E_{4, 6}
    ={1\over 2}\lambda_{8*}( e_{1, 7} - e_{2, 8} + e_{3, 5} - e_{4, 6}),\\
 \kappa_{8*}^+( e_1 e_8)
    &=& E_{1, 8} + E_{2, 7} - E_{3, 6} - E_{4, 5}
    ={1\over 2}\lambda_{8*}( e_{1, 8} + e_{2, 7} - e_{3, 6} - e_{4, 5}),\\
 \kappa_{8*}^+( e_2 e_3)
    &=& E_{1, 4} + E_{2, 3} + E_{5, 8} + E_{6, 7}
    ={1\over 2}\lambda_{8*}( e_{1, 4} + e_{2, 3} + e_{5, 8} + e_{6, 7}),\\
 \kappa_{8*}^+( e_2 e_4)
    &=& -E_{1, 3} + E_{2, 4} + E_{5, 7} - E_{6, 8}
    ={1\over 2}\lambda_{8*}( -e_{1, 3} + e_{2, 4} + e_{5, 7} - e_{6, 8}),\\
 \kappa_{8*}^+( e_2 e_5)
    &=& E_{1, 6} + E_{2, 5} - E_{3, 8} - E_{4, 7}
    ={1\over 2}\lambda_{8*}( e_{1, 6} + e_{2, 5} - e_{3, 8} - e_{4, 7}),\\
 \kappa_{8*}^+( e_2 e_6)
    &=& -E_{1, 5} + E_{2, 6} - E_{3, 7} + E_{4, 8}
    ={1\over 2}\lambda_{8*}( -e_{1, 5} + e_{2, 6} - e_{3, 7} + e_{4, 8}),\\
 \kappa_{8*}^+( e_2 e_7)
   &=& E_{1, 8} + E_{2, 7} + E_{3, 6} + E_{4, 5}
   ={1\over 2}\lambda_{8*}( e_{1, 8} + e_{2, 7} + e_{3, 6} + e_{4, 5}),\\
 \kappa_{8*}^+( e_2 e_8)
   &=& -E_{1, 7} + E_{2, 8} + E_{3, 5} - E_{4, 6}
   ={1\over 2}\lambda_{8*}( -e_{1, 7} + e_{2, 8} + e_{3, 5} - e_{4, 6}),\\
 \kappa_{8*}^+( e_3 e_4)
    &=& E_{1, 2} + E_{3, 4} - E_{5, 6} - E_{7, 8}
    ={1\over 2}\lambda_{8*}( e_{1, 2} + e_{3, 4} - e_{5, 6} - e_{7, 8}),\\
 \kappa_{8*}^+( e_3 e_5)
    &=& E_{1, 7} + E_{2, 8} + E_{3, 5} + E_{4, 6}
    ={1\over 2}\lambda_{8*}( e_{1, 7} + e_{2, 8} + e_{3, 5} + e_{4, 6}),\\
 \kappa_{8*}^+( e_3 e_6)
    &=& -E_{1, 8} + E_{2, 7} + E_{3, 6} - E_{4, 5}
    ={1\over 2}\lambda_{8*}( -e_{1, 8} + e_{2, 7} + e_{3, 6} - e_{4, 5}),\\
 \kappa_{8*}^+( e_3 e_7)
    &=& -E_{1, 5} - E_{2, 6} + E_{3, 7} + E_{4, 8}
    ={1\over 2}\lambda_{8*}( -e_{1, 5} - e_{2, 6} + e_{3, 7} + e_{4, 8}),\\
 \kappa_{8*}^+( e_3 e_8)
    &=& E_{1, 6} - E_{2, 5} + E_{3, 8} - E_{4, 7}
    ={1\over 2}\lambda_{8*}( e_{1, 6} - e_{2, 5} + e_{3, 8} - e_{4, 7}),\\
 \kappa_{8*}^+( e_4 e_5)
    &=& -E_{1, 8} + E_{2, 7} - E_{3, 6} + E_{4, 5}
    ={1\over 2}\lambda_{8*}( -e_{1, 8} + e_{2, 7} - e_{3, 6} + e_{4, 5}),\\
 \kappa_{8*}^+( e_4 e_6)
    &=& -E_{1, 7} - E_{2, 8} + E_{3, 5} + E_{4, 6}
    ={1\over 2}\lambda_{8*}( -e_{1, 7} - e_{2, 8} + e_{3, 5} + e_{4, 6}),\\
 \kappa_{8*}^+( e_4 e_7)
    &=& E_{1, 6} - E_{2, 5} - E_{3, 8} + E_{4, 7}
    ={1\over 2}\lambda_{8*}( e_{1, 6} - e_{2, 5} - e_{3, 8} + e_{4, 7}),\\
 \kappa_{8*}^+( e_4 e_8)
    &=& E_{1, 5} + E_{2, 6} + E_{3, 7} + E_{4, 8}
    ={1\over 2}\lambda_{8*}( e_{1, 5} + e_{2, 6} + e_{3, 7} + e_{4, 8}),\\
 \kappa_{8*}^+( e_5 e_6)
    &=& E_{1, 2} - E_{3, 4} + E_{5, 6} - E_{7, 8}
    ={1\over 2}\lambda_{8*}( e_{1, 2} - e_{3, 4} + e_{5, 6} - e_{7, 8}),\\
 \kappa_{8*}^+( e_5 e_7)
    &=& E_{1, 3} + E_{2, 4} + E_{5, 7} + E_{6, 8}
    ={1\over 2}\lambda_{8*}( e_{1, 3} + e_{2, 4} + e_{5, 7} + e_{6, 8}),\\
 \kappa_{8*}^+( e_5 e_8)
    &=& -E_{1, 4} + E_{2, 3} + E_{5, 8} - E_{6, 7}
    ={1\over 2}\lambda_{8*}( -e_{1, 4} + e_{2, 3} + e_{5, 8} - e_{6, 7}),\\
 \kappa_{8*}^+( e_6 e_7)
    &=& -E_{1, 4} + E_{2, 3} - E_{5, 8} + E_{6, 7}
    ={1\over 2}\lambda_{8*}( -e_{1, 4} + e_{2, 3} - e_{5, 8} + e_{6, 7}),\\
 \kappa_{8*}^+( e_6 e_8)
    &=& -E_{1, 3} - E_{2, 4} + E_{5, 7} + E_{6, 8}
    ={1\over 2}\lambda_{8*}( -e_{1, 3} - e_{2, 4} + e_{5, 7} + e_{6, 8}),\\
 \kappa_{8*}^+( e_7 e_8)
    &=& E_{1, 2} - E_{3, 4} - E_{5, 6} + E_{7, 8}
    ={1\over 2}\lambda_{8*}( e_{1, 2} - e_{3, 4} - e_{5, 6} + e_{7, 8}),
\end{eqnarray*}
This means, in terms of the second diagram in \rf{eq: diagrams Lie algebras},
\begin{eqnarray*}
 \tau_*( e_1 e_2)&=& {1\over 2} (e_1e_2 + e_3e_4 + e_5e_6 + e_7e_8),\\
 \tau_*( e_1 e_3)&=& {1\over 2} (e_1e_3 - e_2e_4 + e_5e_7 - e_6e_8),\\
 \tau_*( e_1 e_4)&=& {1\over 2} (e_1e_4 + e_2e_3 - e_5e_8 - e_6e_7),\\
 \tau_*( e_1 e_5)&=& {1\over 2} (e_1e_5 - e_2e_6 - e_3e_7 + e_4e_8),\\
 \tau_*( e_1 e_6)&=& {1\over 2} (e_1e_6 + e_2e_5 + e_3e_8 + e_4e_7),\\
 \tau_*( e_1 e_7)&=& {1\over 2} (e_1e_7 - e_2e_8 + e_3e_5 - e_4e_6),\\
 \tau_*( e_1 e_8)&=& {1\over 2} (e_1e_8 + e_2e_7 - e_3e_6 - e_4e_5),\\
 \tau_*( e_2 e_3)&=& {1\over 2} (e_1e_4 + e_2e_3 + e_5e_8 + e_6e_7),\\
 \tau_*( e_2 e_4)&=& {1\over 2} (-e_1e_3 + e_2e_4 + e_5e_7 - e_6e_8),\\
 \tau_*( e_2 e_5)&=& {1\over 2} (e_1e_6 + e_2e_5 - e_3e_8 - e_4e_7),\\
 \tau_*( e_2 e_6)&=& {1\over 2} (-e_1e_5 + e_2e_6 - e_3e_7 + e_4e_8),\\
 \tau_*( e_2 e_7)&=& {1\over 2} (e_1e_8 + e_2e_7 + e_3e_6 + e_4e_5),\\
 \tau_*( e_2 e_8)&=& {1\over 2} (-e_1e_7 + e_2e_8 + e_3e_5 - e_4e_6),\\
 \tau_*( e_3 e_4)&=& {1\over 2} (e_1e_2 + e_3e_4 - e_5e_6 - e_7e_8),\\
 \tau_*( e_3 e_5)&=& {1\over 2} (e_1e_7 + e_2e_8 + e_3e_5 + e_4e_6),\\
 \tau_*( e_3 e_6)&=& {1\over 2} (-e_1e_8 + e_2e_7 + e_3e_6 - e_4e_5),\\
 \tau_*( e_3 e_7)&=& {1\over 2} (-e_1e_5 - e_2e_6 + e_3e_7 + e_4e_8),\\
 \tau_*( e_3 e_8)&=& {1\over 2} (e_1e_6 - e_2e_5 + e_3e_8 - e_4e_7),\\
 \tau_*( e_4 e_5)&=& {1\over 2} (-e_1e_8 + e_2e_7 - e_3e_6 + e_4e_5),\\
 \tau_*( e_4 e_6)&=& {1\over 2} (-e_1e_7 - e_2e_8 + e_3e_5 + e_4e_6),\\
 \tau_*( e_4 e_7)&=& {1\over 2} (e_1e_6 - e_2e_5 - e_3e_8 + e_4e_7),\\
 \tau_*( e_4 e_8)&=& {1\over 2} (e_1e_5 + e_2e_6 + e_3e_7 + e_4e_8),\\
 \tau_*( e_5 e_6)&=& {1\over 2} (e_1e_2 - e_3e_4 + e_5e_6 - e_7e_8),\\
 \tau_*( e_5 e_7)&=& {1\over 2} (e_1e_3 + e_2e_4 + e_5e_7 + e_6e_8),\\
 \tau_*( e_5 e_8)&=& {1\over 2} (-e_1e_4 + e_2e_3 + e_5e_8 - e_6e_7),\\
 \tau_*( e_6 e_7)&=& {1\over 2} (-e_1e_4 + e_2e_3 - e_5e_8 + e_6e_7),\\
 \tau_*( e_6 e_8)&=& {1\over 2} (-e_1e_3 - e_2e_4 + e_5e_7 + e_6e_8),\\
 \tau_*( e_7 e_8)&=& {1\over 2} (e_1e_2 - e_3e_4 - e_5e_6 + e_7e_8).
\end{eqnarray*}
In other words, we have defined $\tau_*$ in such a way that
\[\lambda_{8*}\circ \tau_* = \kappa_{8*}^+.\]
As before, using the stardard ordered basis of $\mathfrak{spin}(8)$, we have the matrix representation
{\tiny
\[
\tau_*={1\over 2}
\left(
\begin{array}{cccccccccccccccccccccccccccc}
     1  & 0 & 0 & 0 & 0 & 0 & 0 & 0 & 0 & 0 & 0 & 0 & 0 & 1  & 0 &
    0 & 0 & 0 & 0 & 0 & 0 & 0 & 1  & 0 & 0 & 0 & 0 & 1 \\
     0 & 1  & 0 & 0 & 0 & 0 & 0 & 0 & -1  & 0 & 0 & 0 & 0 & 0 & 0 &
    0 & 0 & 0 & 0 & 0 & 0 & 0 & 0 & 1  & 0 & 0 & -1  & 0\\
     0 & 0 & 1  & 0 & 0 & 0 & 0 & 1  & 0 & 0 & 0 & 0 & 0 & 0 & 0 &
    0 & 0 & 0 & 0 & 0 & 0 & 0 & 0 & 0 & -1  & -1  & 0 & 0\\
     0 & 0 & 0 & 1  & 0 & 0 & 0 & 0 & 0 & 0 & -1  & 0 & 0 & 0 & 0 &
    0 & -1  & 0 & 0 & 0 & 0 & 1  & 0 & 0 & 0 & 0 & 0 & 0\\
     0 & 0 & 0 & 0 & 1  & 0 & 0 & 0 & 0 & 1  & 0 & 0 & 0 & 0 & 0 &
    0 & 0 & 1  & 0 & 0 & 1  & 0 & 0 & 0 & 0 & 0 & 0 & 0\\
     0 & 0 & 0 & 0 & 0 & 1  & 0 & 0 & 0 & 0 & 0 & 0 & -1  & 0 & 1
    & 0 & 0 & 0 & 0 & -1  & 0 & 0 & 0 & 0 & 0 & 0 & 0 & 0\\
     0 & 0 & 0 & 0 & 0 & 0 & 1  & 0 & 0 & 0 & 0 & 1  & 0 & 0 & 0 &
    -1  & 0 & 0 & -1  & 0 & 0 & 0 & 0 & 0 & 0 & 0 & 0 & 0\\
     0 & 0 & 1  & 0 & 0 & 0 & 0 & 1  & 0 & 0 & 0 & 0 & 0 & 0 & 0 &
    0 & 0 & 0 & 0 & 0 & 0 & 0 & 0 & 0 & 1  & 1  & 0 & 0\\
     0 & -1  & 0 & 0 & 0 & 0 & 0 & 0 & 1  & 0 & 0 & 0 & 0 & 0 & 0 &
    0 & 0 & 0 & 0 & 0 & 0 & 0 & 0 & 1  & 0 & 0 & -1  & 0\\
     0 & 0 & 0 & 0 & 1  & 0 & 0 & 0 & 0 & 1  & 0 & 0 & 0 & 0 & 0 &
    0 & 0 & -1  & 0 & 0 & -1  & 0 & 0 & 0 & 0 & 0 & 0 & 0\\
     0 & 0 & 0 & -1  & 0 & 0 & 0 & 0 & 0 & 0 & 1  & 0 & 0 & 0 & 0 &
    0 & -1  & 0 & 0 & 0 & 0 & 1  & 0 & 0 & 0 & 0 & 0 & 0\\
     0 & 0 & 0 & 0 & 0 & 0 & 1  & 0 & 0 & 0 & 0 & 1  & 0 & 0 & 0 &
    1  & 0 & 0 & 1  & 0 & 0 & 0 & 0 & 0 & 0 & 0 & 0 & 0\\
     0 & 0 & 0 & 0 & 0 & -1  & 0 & 0 & 0 & 0 & 0 & 0 & 1  & 0 & 1
    & 0 & 0 & 0 & 0 & -1  & 0 & 0 & 0 & 0 & 0 & 0 & 0 & 0\\
     1  & 0 & 0 & 0 & 0 & 0 & 0 & 0 & 0 & 0 & 0 & 0 & 0 & 1  & 0 &
    0 & 0 & 0 & 0 & 0 & 0 & 0 & -1  & 0 & 0 & 0 & 0 & -1 \\
     0 & 0 & 0 & 0 & 0 & 1  & 0 & 0 & 0 & 0 & 0 & 0 & 1  & 0 & 1
    & 0 & 0 & 0 & 0 & 1  & 0 & 0 & 0 & 0 & 0 & 0 & 0 & 0\\
     0 & 0 & 0 & 0 & 0 & 0 & -1  & 0 & 0 & 0 & 0 & 1  & 0 & 0 & 0 &
    1  & 0 & 0 & -1  & 0 & 0 & 0 & 0 & 0 & 0 & 0 & 0 & 0\\
     0 & 0 & 0 & -1  & 0 & 0 & 0 & 0 & 0 & 0 & -1  & 0 & 0 & 0 & 0
    & 0 & 1  & 0 & 0 & 0 & 0 & 1  & 0 & 0 & 0 & 0 & 0 & 0\\
     0 & 0 & 0 & 0 & 1  & 0 & 0 & 0 & 0 & -1  & 0 & 0 & 0 & 0 & 0 &
    0 & 0 & 1  & 0 & 0 & -1  & 0 & 0 & 0 & 0 & 0 & 0 & 0\\
     0 & 0 & 0 & 0 & 0 & 0 & -1  & 0 & 0 & 0 & 0 & 1  & 0 & 0 & 0 &
    -1  & 0 & 0 & 1  & 0 & 0 & 0 & 0 & 0 & 0 & 0 & 0 & 0\\
     0 & 0 & 0 & 0 & 0 & -1  & 0 & 0 & 0 & 0 & 0 & 0 & -1  & 0 &
    1  & 0 & 0 & 0 & 0 & 1  & 0 & 0 & 0 & 0 & 0 & 0 & 0 & 0\\
     0 & 0 & 0 & 0 & 1  & 0 & 0 & 0 & 0 & -1  & 0 & 0 & 0 & 0 & 0 &
    0 & 0 & -1  & 0 & 0 & 1  & 0 & 0 & 0 & 0 & 0 & 0 & 0\\
     0 & 0 & 0 & 1  & 0 & 0 & 0 & 0 & 0 & 0 & 1  & 0 & 0 & 0 & 0 &
    0 & 1  & 0 & 0 & 0 & 0 & 1  & 0 & 0 & 0 & 0 & 0 & 0\\
     1  & 0 & 0 & 0 & 0 & 0 & 0 & 0 & 0 & 0 & 0 & 0 & 0 & -1  & 0 &
    0 & 0 & 0 & 0 & 0 & 0 & 0 & 1  & 0 & 0 & 0 & 0 & -1 \\
     0 & 1  & 0 & 0 & 0 & 0 & 0 & 0 & 1  & 0 & 0 & 0 & 0 & 0 & 0 &
    0 & 0 & 0 & 0 & 0 & 0 & 0 & 0 & 1  & 0 & 0 & 1  & 0\\
     0 & 0 & -1  & 0 & 0 & 0 & 0 & 1  & 0 & 0 & 0 & 0 & 0 & 0 & 0 &
    0 & 0 & 0 & 0 & 0 & 0 & 0 & 0 & 0 & 1  & -1  & 0 & 0\\
     0 & 0 & -1  & 0 & 0 & 0 & 0 & 1  & 0 & 0 & 0 & 0 & 0 & 0 & 0 &
    0 & 0 & 0 & 0 & 0 & 0 & 0 & 0 & 0 & -1  & 1  & 0 & 0\\
     0 & -1  & 0 & 0 & 0 & 0 & 0 & 0 & -1  & 0 & 0 & 0 & 0 & 0 & 0
    & 0 & 0 & 0 & 0 & 0 & 0 & 0 & 0 & 1  & 0 & 0 & 1  & 0\\
     1  & 0 & 0 & 0 & 0 & 0 & 0 & 0 & 0 & 0 & 0 & 0 & 0 & -1  & 0 &
    0 & 0 & 0 & 0 & 0 & 0 & 0 & -1  & 0 & 0 & 0 & 0 & 1 \\
\end{array}
\right),
\]
}
which one can verify is of order 2.
The map $\tau_*$ has eigenvalues
\[1, -1\]
with multiplicities 21 and 7 respectively.
The eigenspace corresponding to $1$ is generated by
\begin{eqnarray*}
&& \{                    e_{1} e_{4} + e_{2} e_{3},
                     -e_{1} e_{3} + e_{2} e_{4},
                     e_{1} e_{6} + e_{2} e_{5},
                     -e_{1} e_{5} + e_{2} e_{6},
                     e_{1} e_{8} + e_{2} e_{7},
                     -e_{1} e_{7} + e_{2} e_{8},
                     e_{1} e_{2} + e_{3} e_{4},\\
&&                     e_{1} e_{7} + e_{3} e_{5},
                     -e_{1} e_{8} + e_{3} e_{6},
                     -e_{1} e_{5} + e_{3} e_{7},
                     e_{1} e_{6} + e_{3} e_{8},
                     -e_{1} e_{8} + e_{4} e_{5},
                     -e_{1} e_{7} + e_{4} e_{6},
                     e_{1} e_{6} + e_{4} e_{7},\\
&&                     e_{1} e_{5} + e_{4} e_{8},
                     e_{1} e_{2} + e_{5} e_{6},
                     e_{1} e_{3} + e_{5} e_{7},
                     -e_{1} e_{4} + e_{5} e_{8},
                     -e_{1} e_{4} + e_{6} e_{7},
                     -e_{1} e_{3} + e_{6} e_{8},
                     e_{1} e_{2} + e_{7} e_{8}\}
\end{eqnarray*}
which generates a copy of $\mathfrak{spin}(7)\subset\mathfrak{spin}(8)\subset Cl_8^0$, i.e.
\[\{\mbox{+1 eigenspace of $\tau_*$}\}\cong \mathfrak{spin}(7).\]
By taking appropriate sums of these generators we can find
the set of generators \rf{eq: generators g2} of our copy of $\mathfrak{g}_2$.
Moreover, $\mathfrak{g}_2$ is the intersection of the two copies of $\mathfrak{spin}(7)$, i.e.
\[\mathfrak{g}_2= \mathfrak{spin}({\rm span}\{e_2,...,e_8\})\cap \{\mbox{+1 eigenspace of $\tau_*$}\}.\]
One can easlily compute brackets (in Clifford product) of the pair of Lie algebras
$(\mathfrak{spin}(7),\mathfrak{g}_2)$ and check they form a symmetric pair. Since the orbit space
$Spin(7)\cdot e_1=Spin(7)/G_2$ is $7$-dimensional and a submanifold of the $7$-dimensional sphere $S^7=Spin(8)/Spin(7)$, 
we have the classical result \cite{Lawson}
\[{Spin(7)\over G_2}=S^7 .\]

The $7$-dimensional eigenspace corresponding to $-1$ is generated by
\begin{eqnarray*}
&& \{        e_{1} e_{3} + e_{2} e_{4} - e_{5} e_{7} + e_{6} e_{8},
         e_{1} e_{4} - e_{2} e_{3} + e_{5} e_{8} + e_{6} e_{7},
         e_{1} e_{7} + e_{2} e_{8} - e_{3} e_{5} + e_{4} e_{6},
         e_{1} e_{8} - e_{2} e_{7} + e_{3} e_{6} + e_{4} e_{5},\\
&&         -e_{1} e_{6} + e_{2} e_{5} + e_{3} e_{8} + e_{4} e_{7},
         -e_{1} e_{2} + e_{3} e_{4} + e_{5} e_{6} + e_{7} e_{8},
         -e_{1} e_{5} - e_{2} e_{6} - e_{3} e_{7} + e_{4} e_{8}\}.
\end{eqnarray*}

{\bf Remark}.
Note that, by using the bases
\[\{e_1e_2,e_1e_3,\ldots,e_7e_8\}\subset \mathfrak{spin}(8)\]
and  
\[\{E_{1,2}, E_{1,3},\ldots,E_{7,8}\}\subset \mathfrak{so}(8),\]
the matrices representing 
\begin{eqnarray*}
\kappa_{8*}^-:\mathfrak{spin}(8)&\longrightarrow& \mathfrak{so}(8),\\ 
\kappa_{8*}^+:\mathfrak{spin}(8)&\longrightarrow& \mathfrak{so}(8),
\end{eqnarray*} 
with respect to these bases equal the matrices representing $2\sigma_*$ and $2\tau_*$ respectively.
In this way, triality becomes somewhat tautological. 

\subsubsection{Group generated by $\sigma_*$ and $\tau_*$}

\begin{corol}
The endomorphisms $\sigma_*$ and $\tau_*$
generate a copy of the permutation group $S_3$ of three symbols.
\end{corol}
{\em Proof}.
The endomorphisms $\sigma_*$ and $\tau_*$ satisfy
\begin{eqnarray*}
\tau_*^2&=&{\rm Id}_{\mathfrak{spin}(8)},\\
\sigma_*^3&=&{\rm Id}_{\mathfrak{spin}(8)},\\
\sigma_*\tau_*&=&\tau_*\sigma_*^2,\\
\sigma_*^2\tau_*&=&\tau_*\sigma_*,
\end{eqnarray*}
which proves the claim.\qd.

\begin{corol}
The compositions
$\tau_*\sigma_*$, $\sigma_*\tau_*$
are also involutions.
\end{corol}
{\em Proof}. Consider,
\begin{eqnarray*}
(\tau_*\sigma_*)(\tau_*\sigma_*)
    &=&\tau_*(\sigma_*\tau_*)\sigma_*\\
    &=&\tau_*(\tau_*\sigma_*^2)\sigma_*\\
    &=&\tau_*^2\sigma_*^3\\
    &=&{\rm Id}_{\mathfrak{spin}(8)}.
\end{eqnarray*}
\qd

\begin{corol}
The endomorphisms $\lambda_{8*}$, $\kappa_{8*}^+$, $\kappa_{8*}^-$, $\tau_*$ and $\sigma_*$ satisfy
\begin{eqnarray*}
\lambda_{8*}\sigma_*&=& \kappa_{8*}^-,\\
\lambda_{8*}\tau_*&=& \kappa_{8*}^+,\\
\kappa_{8*}^-\tau_*\sigma_*&=& \kappa_{8*}^+,\\
\kappa_{8*}^+\tau_*\sigma_*&=& \kappa_{8*}^-,
\end{eqnarray*}
i.e. the symmetric group $S_3$ generated by $\tau_*$ and $\sigma_*$ permutes the Lie algebra representations $\lambda_{8*}$, $\kappa_{8*}^+$ and $\kappa_{8*}^-$, and the following diagram commutes
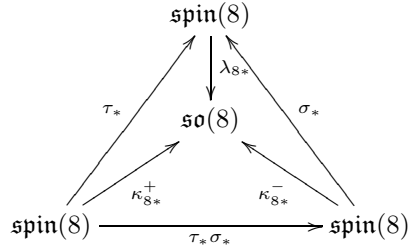
\[
\xymatrix{
& \mathfrak{spin}(8) \ar[d]^{\lambda_{8*}}&\\
& \mathfrak{so}(8)&\\
\mathfrak{spin}(8)\ar[rr]_{\tau_*\sigma_*}\ar[uur]^{\tau_*}\ar[ur]_{\kappa_{8*}^+}&&\mathfrak{spin}(8)
\ar[uul]_{\sigma_*}\ar[ul]^{\kappa_{8*}^-}\\
}
\]

\end{corol}
\qd

We know that the following diagrams also commute
\[
\xymatrix{
 \mathfrak{spin}(8)\ar[r]^{\sigma_*}\ar[d]_{exp} & \mathfrak{spin}(8)\ar[d]^{exp} \\
Spin(8) \ar[r]^{\sigma} & Spin(8)
}
\quad\quad\quad\quad\quad
\xymatrix{
 \mathfrak{spin}(8)\ar[r]^{\tau_*}\ar[d]_{exp} & \mathfrak{spin}(8)\ar[d]^{exp} \\
Spin(8) \ar[r]^{\tau} & Spin(8)
}
\]
\begin{corol}
The symmetric group $S_3$ generated by $\sigma$ and $\tau$ permutes the 
three representations
$\lambda_8$, $\kappa_8^+$ and $\kappa_8^+$
, i.e. the following diagram commutes
\[
\xymatrix{
& Spin(8) \ar[d]^{\lambda_8}&\\
& SO(8)&\\
Spin(8)\ar[uur]^{\tau}\ar[rr]_{\tau\sigma}\ar[ur]_{\kappa_8^+}&&Spin(8)\ar[ul]^{\kappa_8^-}\ar[luu]_{\sigma}\\
}
\]
\end{corol}
\qd

\begin{corol}
We have
\[\{\mbox{$(+1)$-eigenspace of $\tau_*\sigma_*$}\}\cong\mathfrak{spin}({\rm span}(e_2,e_3,e_4,e_5,e_6,e_7,e_8))
=\mathfrak{spin}(7).\]
\end{corol}
{\em Proof}. It is enough to check the effect of $\tau_*\sigma_*$ on the linear generators of $\mathfrak{spin}(8)$.
We have
\begin{eqnarray*}
\tau_*\sigma_*(e_1e_k)&=& -e_1e_k\quad\quad\mbox{for $2\leq k\leq 8$.}\\
\tau_*\sigma_*(e_ie_j)&=& e_ie_j\quad\quad\quad\mbox{for $2\leq i<j\leq 8$.}
\end{eqnarray*}
\qd

\begin{corol}
We have
\[\mathfrak{g}_2=\{\mbox{$(+1)$-eigenspace of $\tau_*$}\}\cap\{\mbox{$(+1)$-eigenspace of $\tau_*\sigma_*^2$}\}.\]
and 
\[\mathfrak{g}_2=\{\mbox{$(+1)$-eigenspace of $\tau_*$}\}\cap\{\mbox{$(+1)$-eigenspace of $\tau_*\sigma_*$}\}.\]
\end{corol}
{\em Proof}.
If $X\in\{\mbox{$(+1)$-eigenspace of $\tau_*$}\}\cap\{\mbox{$(+1)$-eigenspace of $\tau_*\sigma_*^2$}\}$
\begin{eqnarray*}
\tau_*(X)&=& X\\
\tau_*\sigma_*^2(X)&=&X,
\end{eqnarray*}
Then
\begin{eqnarray*}
X&=&\tau_*\sigma_*^2(X)\\
&=&\sigma_*\tau_*(X)\\
&=& \sigma_*(X),
\end{eqnarray*}
which means $X$ is a $(+1)$-eigenvector of $\sigma_*$, thus an element of $\mathfrak{g}_2$.
A dimension count proves the first identity. The second identity is proved similarly.\qd

\begin{corol}
 $\sigma_*^2$ provides an isomorphism between the two copies of $\mathfrak{spin}(7)$. Namely,
\[\sigma_*^2(\{\mbox{$(+1)$-eigenspace of $\tau_*\sigma_*$}\})=
\{\mbox{$(+1)$-eigenspace of $\tau_*$}\}\]
\end{corol}
{\em Proof}.
Let $Y\in \{\mbox{$(+1)$-eigenspace of $\tau_*\sigma_*$}\}$, i.e.
\[\tau_*\sigma_*(Y)=Y.\]
Apply $\sigma_*^2$ to both sides, so that
\[\sigma_*^2\tau_*\sigma_*(Y)= \sigma_*^2(Y),\]
Since
\[\sigma_*^2\tau_*\sigma_*=\sigma_*^2\tau_*\sigma_*^2\sigma_*^2=\sigma_*^2\sigma_*\tau_*\sigma_*^2=\tau_*\sigma_*^2\]
we have
\[\tau_*(\sigma_*^2(Y))= \sigma_*^2(Y).\]
This means that $X=\sigma_*^2(Y)$ is a $(+1)$-eigenvector of $\tau_*$. Since $\tau_*$ is an automorphism, the claim is proved. \qd

\subsubsection{Fundamental $Spin(7)$ $4$-form and $G_2$ $3$-form}
Using the metric, we can dualize the endomorphisms $\kappa_{8*}^+(e_ie_j)$, $2\leq i<j \leq 8$ into 2-forms:
\begin{eqnarray*}
        f_{2, 3} &=&  dx_{1} \wedge dx_{4} + dx_{2} \wedge dx_{3}
           + dx_{5} \wedge dx_{8} + dx_{6} \wedge dx_{7},\\
       f_{2, 4} &=&  -dx_{1} \wedge dx_{3} + dx_{2} \wedge dx_{4}
          + dx_{5} \wedge dx_{7} - dx_{6} \wedge dx_{8},\\
        f_{2, 5} &=&  dx_{1} \wedge dx_{6} + dx_{2} \wedge dx_{5}
           - dx_{3} \wedge dx_{8} - dx_{4} \wedge dx_{7},\\
       f_{2, 6} &=&  -dx_{1} \wedge dx_{5} + dx_{2} \wedge dx_{6}
          - dx_{3} \wedge dx_{7} + dx_{4} \wedge dx_{8},\\
        f_{2, 7} &=&  dx_{1} \wedge dx_{8} + dx_{2} \wedge dx_{7}
           + dx_{3} \wedge dx_{6} + dx_{4} \wedge dx_{5},\\
       f_{2, 8} &=&  -dx_{1} \wedge dx_{7} + dx_{2} \wedge dx_{8}
          + dx_{3} \wedge dx_{5} - dx_{4} \wedge dx_{6},\\
        f_{3, 4} &=&  dx_{1} \wedge dx_{2} + dx_{3} \wedge dx_{4}
           - dx_{5} \wedge dx_{6} - dx_{7} \wedge dx_{8},\\
        f_{3, 5} &=&  dx_{1} \wedge dx_{7} + dx_{2} \wedge dx_{8}
           + dx_{3} \wedge dx_{5} + dx_{4} \wedge dx_{6},\\
       f_{3, 6} &=&  -dx_{1} \wedge dx_{8} + dx_{2} \wedge dx_{7}
          + dx_{3} \wedge dx_{6} - dx_{4} \wedge dx_{5},\\
       f_{3, 7} &=&  -dx_{1} \wedge dx_{5} - dx_{2} \wedge dx_{6}
          + dx_{3} \wedge dx_{7} + dx_{4} \wedge dx_{8},\\
        f_{3, 8} &=&  dx_{1} \wedge dx_{6} - dx_{2} \wedge dx_{5}
           + dx_{3} \wedge dx_{8} - dx_{4} \wedge dx_{7},\\
       f_{4, 5} &=&  -dx_{1} \wedge dx_{8} + dx_{2} \wedge dx_{7}
          - dx_{3} \wedge dx_{6} + dx_{4} \wedge dx_{5},\\
       f_{4, 6} &=&  -dx_{1} \wedge dx_{7} - dx_{2} \wedge dx_{8}
          + dx_{3} \wedge dx_{5} + dx_{4} \wedge dx_{6},\\
        f_{4, 7} &=&  dx_{1} \wedge dx_{6} - dx_{2} \wedge dx_{5}
           - dx_{3} \wedge dx_{8} + dx_{4} \wedge dx_{7},\\
        f_{4, 8} &=&  dx_{1} \wedge dx_{5} + dx_{2} \wedge dx_{6}
           + dx_{3} \wedge dx_{7} + dx_{4} \wedge dx_{8},\\
        f_{5, 6} &=&  dx_{1} \wedge dx_{2} - dx_{3} \wedge dx_{4}
           + dx_{5} \wedge dx_{6} - dx_{7} \wedge dx_{8},\\
        f_{5, 7} &=&  dx_{1} \wedge dx_{3} + dx_{2} \wedge dx_{4}
           + dx_{5} \wedge dx_{7} + dx_{6} \wedge dx_{8},\\
       f_{5, 8} &=&  -dx_{1} \wedge dx_{4} + dx_{2} \wedge dx_{3}
          + dx_{5} \wedge dx_{8} - dx_{6} \wedge dx_{7},\\
       f_{6, 7} &=&  -dx_{1} \wedge dx_{4} + dx_{2} \wedge dx_{3}
          - dx_{5} \wedge dx_{8} + dx_{6} \wedge dx_{7},\\
       f_{6, 8} &=&  -dx_{1} \wedge dx_{3} - dx_{2} \wedge dx_{4}
          + dx_{5} \wedge dx_{7} + dx_{6} \wedge dx_{8},\\
        f_{7, 8} &=&  dx_{1} \wedge dx_{2} - dx_{3} \wedge dx_{4}
           - dx_{5} \wedge dx_{6} + dx_{7} \wedge dx_{8}.
\end{eqnarray*}
We can form the $Spin(7)$-invariant 4-form
\begin{eqnarray*}
\Omega&=&\sum_{2\leq i<j\leq 8}f_{i,j}\wedge f_{i,j}\\
&=&6[-  dx_1  \wedge dx_2 \wedge dx_3 \wedge dx_4+  dx_1 \wedge dx_4 \wedge dx_6 \wedge dx_7
    -  dx_1  \wedge dx_2 \wedge dx_5 \wedge dx_6\\
&&-  dx_1  \wedge dx_3 \wedge dx_5\wedge dx_7-  dx_1  \wedge dx_2 \wedge dx_7 \wedge dx_8
   +  dx_1  \wedge dx_3 \wedge dx_6 \wedge dx_8\\
&&+  dx_1  \wedge dx_4 \wedge dx_5 \wedge dx_8+  dx_2 \wedge dx_4 \wedge dx_5 \wedge dx_7
    -  dx_2 \wedge dx_4 \wedge dx_6 \wedge dx_8\\
&&+  dx_2 \wedge dx_3 \wedge dx_5 \wedge dx_8+  dx_2  \wedge dx_3 \wedge dx_6 \wedge dx_7
   -  dx_3 \wedge dx_4 \wedge dx_5 \wedge dx_6\\
&&-  dx_3 \wedge dx_4 \wedge dx_7 \wedge dx_8 -  dx_5 \wedge dx_6 \wedge dx_7 \wedge dx_8]
\end{eqnarray*}
whose square is a multiple of the 8-dimensional volume form
\[\Omega\wedge\Omega=504\,\, dx_1\wedge dx_2\wedge dx_3\wedge dx_4\wedge dx_5\wedge dx_6\wedge dx_7\wedge dx_8,\]
thus showing that $\Omega$ is non-degenerate.

By integrating out $dx_1$ we get the
$G_2$-invariant 3-form
\begin{eqnarray*}
\varphi
&=&6(
- (dx_{2}\wedge dx_{3}\wedge dx_{4})- (dx_{2}\wedge dx_{5}\wedge dx_{6})- (dx_{2}\wedge dx_{7}\wedge dx_{8})\\
&&- (dx_{3}\wedge dx_{5}\wedge dx_{7})
+ (dx_{4}\wedge dx_{6}\wedge dx_{7})+ (dx_{3}\wedge dx_{6}\wedge dx_{8})+ (dx_{4}\wedge dx_{5}\wedge dx_{8})
).
\end{eqnarray*}

\subsubsection{$\sigma$ and $\tau$ are outer automorphisms}

Now, we will show that $\sigma$ and $\tau$ are outer automorphisms
by showing that they permute the non-trivial central elements of $Spin(8)$, namely
\[-1, {\rm vol}_8, -{\rm vol}_8.\]
Recall that
\[{\rm vol}_8^2=1.\]

Consider, for instance,
\begin{eqnarray*}
\sigma(exp(t e_1e_2))
   &=& \sigma(\cos(t)+\sin(t) e_1e_2)\\
   &=& exp(\sigma_*(te_1e_2))\\
   &=& exp\left({t\over 2} ( -e_1e_2 - e_3e_4 - e_5e_6 - e_7e_8)\right)\\
   &=& exp\left(-{te_1e_2\over 2} \right) exp\left(-{te_3e_4\over 2} \right)
          exp\left(-{te_5e_6\over 2} \right) exp\left(-{te_7e_8\over 2} \right)\\
   &=& (\cos(t/2)-\sin(t/2) e_1e_2)  (\cos(t/2)-\sin(t/2) e_3e_4)
           (\cos(t/2)-\sin(t/2) e_5e_6)  (\cos(t/2)-\sin(t/2) e_7e_8) ,
\end{eqnarray*}
so that
\begin{eqnarray*}
\sigma(e_1e_2)
   &=&\sigma(\cos(\pi/2)+\sin(\pi/2) e_1e_2)\\
   &=& (\cos(\pi/4)-\sin(\pi/4) e_1e_2)  (\cos(\pi/4)-\sin(\pi/4) e_3e_4)
           (\cos(\pi/4)-\sin(\pi/4) e_5e_6)  (\cos(\pi/4)-\sin(\pi/4) e_7e_8) .
\end{eqnarray*}
Note that the calculations are carried out in $Cl_8$ where the exponentials converge.
Similarly, 
\begin{eqnarray*}
\sigma(e_3e_4)
   &=& (\cos(\pi/4)+\sin(\pi/4) e_1e_2)  (\cos(\pi/4)+\sin(\pi/4) e_3e_4)
           (\cos(\pi/4)-\sin(\pi/4) e_5e_6)  (\cos(\pi/4)-\sin(\pi/4) e_7e_8) .
\\
\sigma(e_5e_6)
   &=& (\cos(\pi/4)+\sin(\pi/4) e_1e_2)  (\cos(\pi/4)-\sin(\pi/4) e_3e_4)
           (\cos(\pi/4)+\sin(\pi/4) e_5e_6)  (\cos(\pi/4)-\sin(\pi/4) e_7e_8) .
\\
\sigma(e_7e_8)
   &=& (\cos(\pi/4)+\sin(\pi/4) e_1e_2)  (\cos(\pi/4)-\sin(\pi/4) e_3e_4)
           (\cos(\pi/4)-\sin(\pi/4) e_5e_6)  (\cos(\pi/4)+\sin(\pi/4) e_7e_8) .
\end{eqnarray*}

Now consider
\begin{eqnarray*}
\tau(exp(t e_1e_2))
   &=& \tau(\cos(t)+\sin(t) e_1e_2)\\
   &=& exp(\tau_*(te_1e_2))\\
   &=& exp\left({t\over 2} ( e_1e_2 + e_3e_4 + e_5e_6 + e_7e_8)\right)\\
   &=& exp\left({te_1e_2\over 2} \right) exp\left({te_3e_4\over 2} \right)
          exp\left({te_5e_6\over 2} \right) exp\left({te_7e_8\over 2} \right)\\
   &=& (\cos(t/2)+\sin(t/2) e_1e_2)  (\cos(t/2)+\sin(t/2) e_3e_4)
           (\cos(t/2)+\sin(t/2) e_5e_6)  (\cos(t/2)+\sin(t/2) e_7e_8) ,
\end{eqnarray*}
so that
\begin{eqnarray*}
\tau(e_1e_2)
   &=&\tau(\cos(\pi/2)+\sin(\pi/2) e_1e_2)\\
   &=& (\cos(\pi/4)+\sin(\pi/4) e_1e_2)  (\cos(\pi/4)+\sin(\pi/4) e_3e_4)
           (\cos(\pi/4)+\sin(\pi/4) e_5e_6)  (\cos(\pi/4)+\sin(\pi/4) e_7e_8) .
\end{eqnarray*}
Similarly for all other generators $e_ie_j$ of $\mathfrak{spin}(8)\subset Cl_8^0.$

\begin{corol}
The automorphisms $\tau $ and $\sigma$ are outer automorphisms, of order 2 and 3 respectively,  since they permute the elements of the center of $Spin(8)$, i.e. 
 \begin{eqnarray*}
\sigma(-1)
   &=& {\rm vol}_8,\\
\sigma({\rm vol}_8)
   &=& -{\rm vol}_8,\\
\sigma(-{\rm vol}_8)
   &=& -1,\\
\tau(-1)
   &=& {\rm vol}_8,\\
\tau({\rm vol}_8) &=& -1.
\end{eqnarray*}
\end{corol}
{\em Proof}.
Consider
{\footnotesize
\begin{eqnarray*}
\sigma(-1)
   &=& \sigma(e_1e_2e_1e_2)\\
   &=& \sigma(e_1e_2)\sigma(e_1e_2)\\
   &=& [(\cos(\pi/4)-\sin(\pi/4) e_1e_2)  (\cos(\pi/4)-\sin(\pi/4) e_3e_4) \\
     &&      (\cos(\pi/4)-\sin(\pi/4) e_5e_6)  (\cos(\pi/4)-\sin(\pi/4) e_7e_8)]^2 \\
   &=& (\cos(\pi/4)-\sin(\pi/4) e_1e_2)^2  (\cos(\pi/4)-\sin(\pi/4) e_3e_4)^2 \\
     &&      (\cos(\pi/4)-\sin(\pi/4) e_5e_6)^2  (\cos(\pi/4)-\sin(\pi/4) e_7e_8)^2 \\
   &=& (\cos^2(\pi/4)-\sin^2(\pi/4)-2\sin(\pi_4)\cos(\pi/4) e_1e_2)
          (\cos^2(\pi/4)-\sin^2(\pi/4)-2\sin(\pi_4)\cos(\pi/4) e_3e_4) \\
     &&   (\cos^2(\pi/4)-\sin^2(\pi/4)-2\sin(\pi_4)\cos(\pi/4) e_5e_6)
           (\cos^2(\pi/4)-\sin^2(\pi/4)-2\sin(\pi_4)\cos(\pi/4) e_7e_8) \\
   &=& (\cos(\pi/2)-\sin(\pi/2) e_1e_2) (\cos(\pi/2)-\sin(\pi/2) e_3e_4)
         (\cos(\pi/2)-\sin(\pi/2) e_5e_6) (\cos(\pi/2)-\sin(\pi/2) e_7e_8)    \\
   &=& e_1e_2e_3e_4e_5e_6e_7e_8,\\
\sigma({\rm vol}_8)
   &=& \sigma(e_1e_2)\sigma(e_3e_4)\sigma(e_5e_6)\sigma(e_7e_8)\\
   &=& (\cos(\pi/4)-\sin(\pi/4) e_1e_2)  (\cos(\pi/4)-\sin(\pi/4) e_3e_4)
           (\cos(\pi/4)-\sin(\pi/4) e_5e_6)  (\cos(\pi/4)-\sin(\pi/4) e_7e_8) \\
     &&(\cos(\pi/4)+\sin(\pi/4) e_1e_2)  (\cos(\pi/4)+\sin(\pi/4) e_3e_4)
           (\cos(\pi/4)-\sin(\pi/4) e_5e_6)  (\cos(\pi/4)-\sin(\pi/4) e_7e_8)\\
     && (\cos(\pi/4)+\sin(\pi/4) e_1e_2)  (\cos(\pi/4)-\sin(\pi/4) e_3e_4)
           (\cos(\pi/4)+\sin(\pi/4) e_5e_6)  (\cos(\pi/4)-\sin(\pi/4) e_7e_8)\\
     && (\cos(\pi/4)+\sin(\pi/4) e_1e_2)  (\cos(\pi/4)-\sin(\pi/4) e_3e_4)
           (\cos(\pi/4)-\sin(\pi/4) e_5e_6)  (\cos(\pi/4)+\sin(\pi/4) e_7e_8)\\
   &=&    (\cos(\pi/4)+\sin(\pi/4) e_1e_2)^2  (\cos(\pi/4)-\sin(\pi/4) e_3e_4)^2
             (\cos(\pi/4)-\sin(\pi/4) e_5e_6)^2 (\cos(\pi/4)-\sin(\pi/4) e_7e_8)^2\\
   &=& (\cos(\pi/2) + \sin(pi/2) e_1e_2) (\cos(\pi/2) - \sin(pi/2) e_1e_2)
           (\cos(\pi/2) - \sin(pi/2) e_1e_2) (\cos(\pi/2) - \sin(pi/2) e_1e_2)\\
   &=& -e_1e_2e_3e_4e_5e_6e_7e_8,\\
\sigma(-{\rm vol}_8)
   &=& \sigma(-1)\sigma(e_1e_2e_3e_4e_5e_6e_7e_8)\\
   &=& (e_1e_2e_3e_4e_5e_6e_7e_8)(-e_1e_2e_3e_4e_5e_6e_7e_8)\\
   &=& -1.
\end{eqnarray*}

On the other hand, we also have
\begin{eqnarray*}
\tau(-1)
   &=& \tau(e_1e_2e_1e_2)\\
   &=& \tau(e_1e_2)\tau(e_1e_2)\\
   &=& [(\cos(\pi/4)+\sin(\pi/4) e_1e_2)  (\cos(\pi/4)+\sin(\pi/4) e_3e_4) \\
     &&      (\cos(\pi/4)+\sin(\pi/4) e_5e_6)  (\cos(\pi/4)+\sin(\pi/4) e_7e_8)]^2 \\
   &=& (\cos(\pi/4)+\sin(\pi/4) e_1e_2)^2  (\cos(\pi/4)+\sin(\pi/4) e_3e_4)^2 \\
     &&      (\cos(\pi/4)+\sin(\pi/4) e_5e_6)^2  (\cos(\pi/4)+\sin(\pi/4) e_7e_8)^2 \\
   &=& (\cos^2(\pi/4)-\sin^2(\pi/4)+2\sin(\pi_4)\cos(\pi/4) e_1e_2)
          (\cos^2(\pi/4)-\sin^2(\pi/4)+2\sin(\pi_4)\cos(\pi/4) e_3e_4) \\
     &&   (\cos^2(\pi/4)-\sin^2(\pi/4)+2\sin(\pi_4)\cos(\pi/4) e_5e_6)
           (\cos^2(\pi/4)-\sin^2(\pi/4)+2\sin(\pi_4)\cos(\pi/4) e_7e_8) \\
   &=& (\cos(\pi/2)+\sin(\pi/2) e_1e_2) (\cos(\pi/2)+\sin(\pi/2) e_3e_4)
         (\cos(\pi/2)+\sin(\pi/2) e_5e_6) (\cos(\pi/2)+\sin(\pi/2) e_7e_8)    \\
   &=& e_1e_2e_3e_4e_5e_6e_7e_8,\\
\tau({\rm vol}_8) &=& -1.
\end{eqnarray*}
}
\qd

\subsubsection{Octonions}

The relationship between $Spin(8)$ and the Octonions is well known (see \cite{Baez,Lounesto}).
In this subsection, as an example of the use of the binary encoding, we recover a multiplication table of the normed division 
algebra of Octonions using the $Spin(8)$ representations. The idea is to consider Clifford multiplication 
and the three real representations of $Spin(8)$ at the same time. We follow \cite{Baez} and consider 
Clifford multiplication as a trilinear map (a "triality" as defined by Adams in \cite{Adams2}) 
and dualization to get a bilinear map $\mathbb R^8\times \tilde{\Delta}_8^+\rightarrow \tilde{\Delta}_8^-$. 
By identifying the three spaces with a single space (in a suitable way) one can define a product on it using this bilinear map.

Consider the basis of positive real spinors given by
\[\beta^+=\{u_0-u_{15},iu_0+iu_{15},u_3+u_{12},-iu_3+iu_{12},-u_5+u_{10},iu_5+iu_{10},u_6+u_9,iu_6-iu_9\}
\subset \tilde{\Delta}_8^+\]
Let us consider Clifford multiplication as a bilinear map
\[\mathbb R^8\times \tilde{\Delta}_8^+\rightarrow \tilde{\Delta}_8^-.\]
For this subsection, let us denote $\{v_0,\dots,v_7\}$ the standard ordered basis of $\mathbb R^8$.
The Clifford multiplication table is the following:
\[
\begin{array}{|r|rrrrrrrr|}
\hline
& u_0-u_{15}&iu_0+iu_{15}&u_3+u_{12}&-iu_3+iu_{12}&-u_5+u_{10}&iu_5+iu_{10}&u_6+u_9&iu_6-iu_9\\
\hline
v_0& iu_{1}-iu_{14} & -u_{1}-u_{14} & iu_{2}+iu_{13} & -u_{2}+u_{13} & iu_{4}-iu_{11} & -u_{4}-u_{11} & iu_{7}+iu_{8} & -u_{7}+u_{8}\\
v_1&u_{1}+u_{14} & iu_{1}-iu_{14} & -u_{2}+u_{13} & -iu_{2}-iu_{13} & -u_{4}-u_{11} & -iu_{4}+iu_{11} & u_{7}-u_{8} & iu_{7}+iu_{8} \\
v_2&-iu_{2}-iu_{13} & u_{2}-u_{13} & iu_{1}-iu_{14} & -u_{1}-u_{14} & iu_{7}+iu_{8} & -u_{7}+u_{8} & -iu_{4}+iu_{11} &   u_{4}+u_{11} \\
v_3&-u_{2}+u_{13} & -iu_{2}-iu_{13} & -u_{1}-u_{14} & -iu_{1}+iu_{14} & u_{7}-u_{8} & iu_{7}+iu_{8} & u_{4}+u_{11} &   iu_{4}-iu_{11} \\
v_4&iu_{4}-iu_{11} & -u_{4}-u_{11} & iu_{7}+iu_{8} & -u_{7}+u_{8} & -iu_{1}+iu_{14} & u_{1}+u_{14} & -iu_{2}-iu_{13} & u_{2}-u_{13} \\
v_5&u_{4}+u_{11} & iu_{4}-iu_{11} & u_{7}-u_{8} & iu_{7}+iu_{8} & u_{1}+u_{14} & iu_{1}-iu_{14} & u_{2}-u_{13} &   iu_{2}+iu_{13}  \\
v_6&-iu_{8}-iu_{7} & u_{8}-u_{7} & -iu_{11}+iu_{4} & u_{11}+u_{4} & -iu_{13}-iu_{2} & u_{13}-u_{2} & -iu_{14}+iu_{1}  & u_{14}+u_{1} \\
v_7&-u_{8}+u_{7} & -iu_{8}-iu_{7} & -u_{11}-u_{4} & -iu_{11}+iu_{4} & -u_{13}+u_{2} & -iu_{13}-iu_{2} & -u_{14}-u_{1} & -iu_{14}+iu_{1} \\
\hline
\end{array}
\]
Now let
\begin{eqnarray*}
\beta^-&=&\{
                      iu_{1}-iu_{14},
                       -u_{1}-u_{14},
                       iu_{2}+iu_{13},
                       -u_{2}+ u_{13},
                        iu_{4}-iu_{11},
                       -u_{4}- u_{11},
                        iu_{7}+iu_{8},
                       -u_{7}+u_{8}\}
.
\end{eqnarray*}
By labeling the elements of the ordered bases $\beta^+=\{\psi_0,\dots,\psi_7\}$  and  $\beta^-=\{\phi_0,\dots,\phi_7\}$,
the Clifford multiplication table now reads as follows:
\[
\begin{array}{|r|rrrrrrrr|}
\hline
&  \psi_0 & \psi_1 & \psi_2 & \psi_3 & \psi_4 & \psi_5 & \psi_6 & \psi_7 \\
\hline
v_0&\phi_0 &\phi_1 & \phi_2 & \phi_3 & \phi_4& \phi_5 & \phi_6 & \phi_7\\
v_1&-\phi_1 &\phi_0 & \phi_3 & -\phi_2 & \phi_5 & -\phi_4& -\phi_7 & \phi_6 \\
v_2& -\phi_2 & -\phi_3 &\phi_0 &\phi_1 & \phi_6 & \phi_7 & -\phi_4&   -\phi_5 \\
v_3& \phi_3 & -\phi_2 &\phi_1 & -\phi_0 & -\phi_7 & \phi_6 & -\phi_5 &   \phi_4\\
v_4& \phi_4& \phi_5 & \phi_6 & \phi_7 & -\phi_0 &-\phi_1 & -\phi_2 & -\phi_3 \\
v_5& -\phi_5 & \phi_4& -\phi_7 & \phi_6 &-\phi_1 &\phi_0 & -\phi_3 &  \phi_2 \\
v_6&-\phi_6 &\phi_7 & \phi_4 & -\phi_5 & -\phi_2& \phi_3 & \phi_0 & -\phi_1 \\
v_7&-\phi_7 &-\phi_6 & \phi_5 & \phi_4 & -\phi_3& -\phi_2 & \phi_1 & \phi_0\\
\hline
\end{array}
\]
We can recover the multiplication table of the octonions by identifying $\mathbb R^8$,  $ \tilde{\Delta}_8^+$ and $\tilde{\Delta}_8^-$ with a single vector space $\mathbb O={\rm span}\{\hat {e}_0,\hat {e}_1,\dots, \hat {e}_7\}$ in the following way. We identify $v_0$, $\psi_0$ and $\phi_0$ with the identity $\hat {e}_0$ of $\mathbb O$. We also identify $\psi_i$ with  $\hat {e}_i$. In this way, we have that  $\phi_i=v_0 \psi_i=\hat {e}_0\cdot \hat {e}_i=\hat {e}_i$. We also have that $v_1\psi_0=-\phi_1$ so that $v_1$ should be identify with $-\hat {e}_1$. In the same way
$v_2$, $v_3$,  $v_4$, $v_5$, $v_6$ and $v_7$ should be identify with $-\hat {e}_2$, $\hat {e}_3$  $\hat {e}_4$, $-\hat {e}_5$, $-\hat {e}_6$ and $-\hat {e}_7$. Then the multiplication table (of the octonions) reads as follows
\[\begin{array}{|r|rrrrrrrr|}
\hline
& \hat {e}_0 & \hat {e}_1 & \hat {e}_2 & \hat {e}_3 & \hat {e}_4 & \hat {e}_5 & \hat {e}_6 & \hat {e}_7 \\
\hline
\hat {e}_0&\hat {e}_0 &\hat {e}_1 & \hat {e}_2 & \hat {e}_3 & \hat {e}_4& \hat {e}_5 & \hat {e}_6 & \hat {e}_7\\
\hat {e}_1&\hat {e}_1 &-\hat {e}_0 & -\hat {e}_3 & \hat {e}_2 & -\hat {e}_5 & \hat {e}_4& \hat {e}_7 & -\hat {e}_6 \\
\hat {e}_2& \hat {e}_2 & \hat {e}_3 &-\hat {e}_0 &-\hat {e}_1 & -\hat {e}_6 & -\hat {e}_7 & \hat {e}_4&   \hat {e}_5 \\
\hat {e}_3& \hat {e}_3 & -\hat {e}_2 &\hat {e}_1 & -\hat {e}_0 & -\hat {e}_7 & \hat {e}_6 & -\hat {e}_5 &   \hat {e}_4\\
\hat {e}_4& \hat {e}_4& \hat {e}_5 & \hat {e}_6 & \hat {e}_7 & -\hat {e}_0 &-\hat {e}_1 & -\hat {e}_2 & -\hat {e}_3 \\
\hat {e}_5& \hat {e}_5 & -\hat {e}_4& \hat {e}_7 & -\hat {e}_6 &\hat {e}_1 &-\hat {e}_0 & \hat {e}_3 &  -\hat {e}_2 \\
\hat {e}_6&\hat {e}_6 &-\hat {e}_7 & -\hat {e}_4 & \hat {e}_5 & \hat {e}_2& -\hat {e}_3 & -\hat {e}_0 & \hat {e}_1 \\
\hat {e}_7&\hat {e}_7 &\hat {e}_6 & -\hat {e}_5 & -\hat {e}_4 & \hat {e}_3& \hat {e}_2 & -\hat {e}_1 & -\hat {e}_0\\
\hline
\end{array}
\]

{\bf Remark}. One can actually do the same in the case of $\mathbb{R}^4$, $\tilde\Delta_4^+$ and $\tilde\Delta_4^-$ to
recover the quaternion multiplication table.

\subsection{Vector fields on spheres}

Classical results of Hurwitz, Radon and Adams \cite{Hurwitz, Adams} tell us the maximal number of independent vector fields a
sphere can admit, which is given in terms of the Hurwitz Radon numbers.
In this subsection, we will give explicit expressions, using Clifford algebras and the binary code,
for a maximal set orthogonal linearly independent vector fields on spheres (compare with \cite{Piccinni,Ognikyan} for recent work). We follow the ideas described in \cite{Husemoller,Lawson} using Clifford multiplication, but we will multiply by elements of 
$\mathfrak{spin}(r)$ instead of multiplying by elements of the standard basis of $\mathbb{R}^r$.

The idea is as follows: Let $N\in \mathbb N$ and suppose that $\mathbb{R}^N$ is a non-trivial representation of $Cl_r^0$ (not necessarily irreducible). Since Clifford multiplication by unit vectors is an orthogonal transformation on the space of real spinors, 
for every $Z\in S^{N-1}$ and $2\leq j\leq r$, the $r-1$ vectors $V_{j-1}(Z):=e_1e_jZ$ form an orthonormal set tangent to the sphere 
at $Z$. 
By \cite{Adams}, if $r$ is the maximum integer such that $\mathbb{R}^N$ is a non-trivial representation of $Cl_r^0$, then 
the set of independent vector fields is maximal. 
Although the calculations below may seem cumbersome due to the slightly more complicated form of the basic vectors of the real Spin representations \cite{Arizmendi}, the main point is that our expressions provide a general way to produce explicitly the vector fields. As a
concrete example, we compute the $9$ vector fields on $S^{31}$.

Thus, let us suppose that $Cl_r^0$ is represented
on $\mathbb{R}^N$, for some $N\in\mathbb{N}$,
in such a way that each bivector $e_ie_j$ is mapped to an antisymmetric endomorphism $J_{ij}$ satisfying
\begin{equation}
J_{ij}^2 = -{\rm Id}_{\mathbb{R}^N}.\label{eq:almost-complex-structures}
\end{equation}

\begin{itemize}
\item
If $r\not\equiv 0
\,\,\,({\rm mod}\,\,\,\, 4)$, $r>1$, $\mathbb{R}^N$ decomposes into a sum of irreducible representations of $Cl_r^0$.
Since this algebra is simple, such irreducible representations can only be trivial or
copies of the standard representation
$\tilde\Delta_r$ of $Cl_r^0$ (cf. \cite{Lawson}). Due to
\rf{eq:almost-complex-structures}, there are no trivial summands in such a decomposition so that
\begin{eqnarray*}
\mathbb{R}^N
&=& \underbrace{\tilde\Delta_r\oplus\cdots\oplus \tilde\Delta_r}_{m
\,\,\,times} .
\end{eqnarray*}
By restricting to $\mathfrak{spin}(r)\subset Cl_r^0$,
\[\mathbb{R}^N =\tilde\Delta_r \otimes_{\mathbb{R}}\mathbb{R}^m\]
we see that $\mathfrak{spin}(r)$ has an isomorphic image
\[\widehat{\mathfrak{spin}(r)}=\mathfrak{spin}(r)\otimes \{{\rm Id}_{\mathbb{R}^m}\}\subset
\mathfrak{so}(d_rm),\]
which is a subalgebra of $\mathfrak{so}(d_rm)$.
Note that
\[J_{ij}=[\kappa_r(e_ie_j)\otimes{\rm Id}_{\mathbb{R}^m}]\]
for $1\leq i<j\leq r$.

\item
If, on the other hand, $r\equiv 0 \,\,\,({\rm mod}\,\,\,\, 4)$,
\[\hat\Delta_r = \tilde{\Delta}_r^+ \oplus \tilde{\Delta}_r^-,\]
the sum of two inequivalent irreducible representationsm, and
\begin{eqnarray*}
\mathbb{R}^N
&=& \tilde\Delta_r^+\otimes_{\mathbb{R}} \mathbb{R}^{m_1} \oplus \tilde\Delta_r^-\otimes_{\mathbb{R}} \mathbb{R}^{m_2}.
\end{eqnarray*}
as a representation to $\mathfrak{spin}(r)\subset Cl_r^0$, and
we see that $\mathfrak{spin}(r)$ has an isomorphic image
\[\widehat{\mathfrak{spin}(r)}=\{\kappa_n^+(g)\otimes ({\rm Id}_{\mathbb{R}^{m_1}}\oplus
\mathbf{0}_{m_2\times m_2})\oplus
\kappa_n^-(g)\otimes(\mathbf{0}_{m_1\times m_1}\oplus {\rm Id}_{\mathbb{R}^{m_2}})\,|\,g\in \mathfrak{spin}(r)\}
\subset
\mathfrak{so}(d_rm_1+d_rm_2).\]
Note that
\[J_{ij}=[\kappa_r^+(e_ie_j)\otimes{\rm Id}_{\mathbb{R}^{m_1}}]\oplus [\kappa_r^-(e_ie_j)\otimes{\rm Id}_{\mathbb{R}^{m_2}}]\]
for $1\leq i<j\leq r$.

\end{itemize}

Given a point $Z$ in the sphere of $\tilde\Delta_r\otimes_{\mathbb{R}}\mathbb{R}^m$ or
$\tilde\Delta_r^\pm\otimes_{\mathbb{R}}\mathbb{R}^m$, the corresponding values of the
vector fields at $Z$ will be given by
\[ [\kappa_r(e_1e_2)\otimes{\rm Id}_{m}](Z),...,[\kappa_r(e_1e_r)\otimes{\rm Id}_{m}](Z),\]
or
\[[\kappa_r^\pm(e_1e_2)\otimes{\rm Id}_{m}] (Z),...,[\kappa_r^\pm(e_1e_r)\otimes{\rm Id}_{m}](Z),\]
respectively, where ${\rm Id}_{m}:={\rm Id}_{\mathbb{R}^m}$.

\subsubsection{Calculations in $\Delta_r$}

First recall that
\begin{eqnarray*}
 e_{1}u_b
&=& i u_{b+(-1)^{b_0}}.
\end{eqnarray*}
Now, if
$p=2\leq r$, $j=[{p+1\over 2}]= 1$,
\begin{eqnarray*}
 e_1e_2\cdot u_a
   &=& i(-1)^{a_{0}} u_{a},
\end{eqnarray*}
if
$3\leq p\leq 2[{r\over 2}]$, $j=[{p+1\over 2}]\geq 2$,
\begin{eqnarray*}
 e_1e_p\cdot u_a
   &=& i^{1-p}(-1)^{2j-1+\sum_{s=0}^{j-2}a_s+a_{j-1}(-2j+p+1)} u_{a+(-1)^{a_{j-1}}2^{j-1}+(-1)^{a_0}}
\end{eqnarray*}
and if $p=r=2k+1$,
\begin{eqnarray*}
e_1 e_{r}u_a
&=& (-1)^{[r/2]+1+\sum_{l=0}^{[r/2]-1}a_l}  u_{a+(-1)^{a_0}}  .
\end{eqnarray*}
We also have the following expressions for the real and quaternionic structures:
For $r=0,1,4,5 \,\,(\mod 8)$ and $q=[{r/4}]$
\begin{eqnarray*}
\gamma_r(u_a)
   &=&(-i)^q (-1)^{\sum_{t=1}^qa_{2t-1}} u_{2^{2q}-1-a}.
\end{eqnarray*}
For $r=2,3,6,7 \,\,(\mod 8)$ and $q=[{r/4}]$
\begin{eqnarray*}
\gamma_r(u_a)
   &=&(-i)^{q+1} (-1)^{\sum_{t=0}^qa_{2t}} u_{2^{2q+1}-1-a}.
\end{eqnarray*}

\subsubsection{The vector fields}
Due to the coincidence of dimensions of the real Spin representations
\[d_{8q+3}=d_{8q+4},\]
and
\[d_{8q+5}=d_{8q+6}=d_{8q+7}=d_{8q+8},\]
we only need to consider the cases $r\equiv 0,1,2,4 \,\,(\mod 8)$. 
For example, if $R^N$ is a representation space of $Cl_{8q+3}^0$, then  $R^N$ is also a representation of  
$Cl_{8q+4}^0$ and therefore $8q+3$ is not maximal.
Let us consider bases for the spaces $\tilde\Delta_r$ and $\tilde\Delta_r^\pm$:
\begin{itemize}
\item Case $r\equiv 0\,\,(\mod 8)$: A basis for $\tilde\Delta_r^+$ is given by
\[\left\{u_a+\gamma_r(u_a),iu_a+\gamma_r(iu_a)\,|\, a=0,...,2^{r/2-1}-1, \sum a_l\equiv 0 \,\,(\mod 2)\right\}\]

\item Case $r\equiv 1\,\,(\mod 8)$:  A basis for $\tilde\Delta_r$ is given by
\[\left\{u_a+\gamma_r(u_a),iu_a+\gamma_r(iu_a)\,|\, a=0,...,2^{[r/2]-1}-1\right\}\]

\item Case $r\equiv 2\,\,(\mod 8)$:  A basis for $\tilde\Delta_r$ is given by
\[\left\{u_a,iu_a\,|\, a=0,...,2^{r/2}-1,\sum a_l\equiv 0 \,\,(\mod 2)\right\}\]

\item Case $r\equiv 4\,\,(\mod 8)$:  A basis for $\tilde\Delta_r^+$  is given by 
\[\left\{u_a,iu_a\,|\, a=0,...,2^{r/2}-1,\sum a_l\equiv 0 \,\,(\mod 2)\right\}\]

\end{itemize}

\subsubsection*{Case $r\equiv 0\,\,(\mod 8)$}
For any point $Z\in S(\tilde\Delta_r^+\otimes_{\mathbb{R}}\mathbb{R}^m )=S^{d_r m}$,
\[Z= \sum_{h=1}^m\sum_{\{a=0,...,2^{r/2-1}-1\,:\,\sum a_l\equiv 0 \,\,({\rm mod }2)\}} X_{a,h}(u_a\otimes v_h) + Y_{a,h} (iu_a\otimes v_h), \]
where $\{v_1,...,v_m\}$ is an orthonormal basis of $\mathbb{R}^m$ and $X_{a,h},Y_{a,h}\in\mathbb{R}$,
there are
$r-1$ point-wise linearly independent vector fields given as follows.
For $p=2$,
\begin{eqnarray*}
&&[\kappa_r(e_1e_2)\otimes{\rm Id}_m] \left(\sum_{h=1}^m\sum_{\{a=0,...,2^{r/2-1}-1,\,\sum a_l\equiv 0 \,\,({\rm mod }2)\}}
X_{a,h}(u_a+\gamma_r(u_a))\otimes v_h) + Y_{a,h} (iu_a+\gamma_r(iu_a))\otimes v_h)\right)\\
&=&\sum_{h=1}^m\sum_{\{a=0,...,2^{r/2-1}-1,\,\sum a_l\equiv 0 \,\,({\rm mod }2)\}}
(-1)^{a_{0}}(-Y_{a,h} ( u_{a}+\gamma_r( u_{a}))\otimes v_h)+X_{a,h}(i u_{a}+\gamma_r(i u_{a}))\otimes v_h) ).
\end{eqnarray*}
For $3\leq p\leq r$,
\begin{eqnarray*}
&&[\kappa_r(e_1e_p)\otimes{\rm Id}_m] \left(\sum_{\{a=0,...,2^{r/2-1}-1,\,\sum a_l\equiv 0 \,\,({\rm mod }2)\}}
X_{a,h}(u_a+\gamma_r(u_a))\otimes v_h) + Y_{a,h} (iu_a+\gamma_r(iu_a))\otimes v_h)\right)\\
&=&\sum_{h=1}^m\sum_{\{a=0,...,2^{r/2-1}-1,\,\sum a_l\equiv 0 \,\,({\rm mod }2)\}}
(-1)^{2j-1+\sum_{s=0}^{j-2}a_s+a_{j-1}(-2j+p+1)}\\
&&(X_{a,h}(i^{1-p} u_{a+(-1)^{a_{j-1}}2^{j-1}+(-1)^{a_0}}+\gamma_r(i^{1-p} u_{a+(-1)^{a_{j-1}}2^{j-1}+(-1)^{a_0}}))\otimes v_h) \\
&& -Y_{a,h} (i^{-p} u_{a+(-1)^{a_{j-1}}2^{j-1}+(-1)^{a_0}}+\gamma_r(i^{-p} u_{a+(-1)^{a_{j-1}}2^{j-1}+(-1)^{a_0}}))\otimes v_h)\\
\end{eqnarray*}

\subsubsection*{Case $r\equiv 1\,\,(\mod 8)$}
For any point) $Z\in S(\tilde\Delta_r\otimes_{\mathbb{R}}\mathbb{R}^m )=S^{d_r m}$,
\[Z= \sum_{h=1}^m\sum_{a=0}^{2^{[r/2]-1}-1} X_{a,h}(u_a\otimes v_h) + Y_{a,h} (iu_a\otimes v_h), \]
where $\{v_1,...,v_m\}$ is an orthonormal basis of $\mathbb{R}^m$ and $X_{a,h},Y_{a,h}\in\mathbb{R}$,
there are
$r-1$ point-wise linearly independent vector fields given as follows.
For $p=2$,
\begin{eqnarray*}
&&[\kappa_r(e_1e_2)\otimes{\rm Id}_m] \left(\sum_{h=1}^m\sum_{a=0}^{2^{[r/2]-1}-1}
X_{a,h}(u_a+\gamma_r(u_a))\otimes v_h) + Y_{a,h} (iu_a+\gamma_r(iu_a))\otimes v_h)\right)\\
&=&\sum_{h=1}^m\sum_{a=0}^{2^{[r/2]-1}-1}
(-1)^{a_{0}}(-Y_{a,h} ( u_{a}+\gamma_r( u_{a}))\otimes v_h)+X_{a,h}(i u_{a}+\gamma_r(i u_{a}))\otimes v_h) ).
\end{eqnarray*}
For $3\leq p\leq r-1$,
\begin{eqnarray*}
&&[\kappa_r(e_1e_p)\otimes{\rm Id}_m] \left(\sum_{h=1}^m\sum_{a=0}^{2^{[r/2]-1}-1}
X_{a,h}(u_a+\gamma_r(u_a))\otimes v_h) + Y_{a,h} (iu_a+\gamma_r(iu_a))\otimes v_h)\right)\\
&=&\sum_{h=1}^m\sum_{a=0}^{2^{[r/2]-1}-1}
(-1)^{2j-1+\sum_{s=0}^{j-2}a_s+a_{j-1}(-2j+p+1)}(X_{a,h}(i^{1-p} u_{a+(-1)^{a_{j-1}}2^{j-1}+(-1)^{a_0}}+\gamma_r(i^{1-p} u_{a+(-1)^{a_{j-1}}2^{j-1}+(-1)^{a_0}}))\otimes v_h) \\
&& -Y_{a,h} (i^{-p} u_{a+(-1)^{a_{j-1}}2^{j-1}+(-1)^{a_0}}+\gamma_r(i^{-p} u_{a+(-1)^{a_{j-1}}2^{j-1}+(-1)^{a_0}}))\otimes v_h)\\
\end{eqnarray*}
For $p=r$,
\begin{eqnarray*}
&&[\kappa_r(e_1e_r)\otimes{\rm Id}_m] \left(\sum_{h=1}^m\sum_{a=0}^{2^{[r/2]-1}-1}
X_{a,h}(u_a+\gamma_r(u_a))\otimes v_h) + Y_{a,h} (iu_a+\gamma_r(iu_a))\otimes v_h)\right)\\
&=& \sum_{h=1}^m\sum_{a=0}^{2^{[r/2]-1}-1}
(-1)^{[r/2]+1+\sum_{l=0}^{[r/2]-1}a_l}(X_{a,h}(  u_{a+(-1)^{a_0}}+\gamma_r(  u_{a+(-1)^{a_0}}))\otimes v_h) \\
&&+ Y_{a,h} (i  u_{a+(-1)^{a_0}}+\gamma_r(i  u_{a+(-1)^{a_0}}))\otimes v_h).
\end{eqnarray*}

\subsubsection*{Case $r\equiv 2,4\,\,(\mod 8)$}

For any point $Z\in S(\tilde\Delta_r\otimes_{\mathbb{R}}\mathbb{R}^m )=S^{d_r m}$ if $r\equiv 2\,\,(\mod 8)$
(resp. $Z\in S(\tilde\Delta_r^+\otimes_{\mathbb{R}}\mathbb{R}^m )=S^{d_r m}$ if $r\equiv 4\,\,(\mod 8)$),
\[Z= \sum_{h=1}^m\sum_{\{a=0,...,2^{r/2}-1,\,\sum a_l\equiv 0 \,\,({\rm mod }2)\}} X_{a,h}(u_a\otimes v_h) + Y_{a,h} (iu_a\otimes v_h), \]
where $\{v_1,...,v_m\}$ is an orthonormal basis of $\mathbb{R}^m$ and $X_{a,h},Y_{a,h}\in\mathbb{R}$,
there are
$r-1$ point-wise linearly independent vector fields given as follows.
For $p=2$,
\begin{eqnarray*}
&&[\kappa_r(e_1e_2)\otimes{\rm Id}_m] \left(\sum_{h=1}^m\sum_{\{a=0,...,2^{r/2}-1,\,\sum a_l\equiv 0 \,\,({\rm mod }2)\}} X_{a,h}(u_a\otimes v_h) + Y_{a,h} (iu_a\otimes v_h)\right)
\\
   &=& \sum_{h=1}^m \sum_{\{a=0,...,2^{r/2}-1,\,\sum a_l\equiv 0 \,\,({\rm mod }2)\}}  (-1)^{a_{0}}(-Y_{a,h}  (u_{a}\otimes v_h) + X_{a,h} (iu_{a}\otimes v_h)).
\end{eqnarray*}
For $3\leq p\leq r$,
\begin{eqnarray*}
&&[\kappa_r(e_1e_p)\otimes{\rm Id}_m] \left(\sum_{h=1}^m\sum_{\{a=0,...,2^{r/2}-1,\,\sum a_l\equiv 0 \,\,({\rm mod }2)\}} X_{a,h}(u_a\otimes v_h) + Y_{a,h} (iu_a\otimes v_h)\right)\\
   &=& \sum_{h=1}^m\sum_{\{a=0,...,2^{r/2}-1,\,\sum a_l\equiv 0 \,\,({\rm mod }2)\}} (-1)^{2j-1+\sum_{s=0}^{j-2}a_s+a_{j-1}(-2j+p+1)}(X_{a,h}(i^{1-p} u_{a+(-1)^{a_{j-1}}2^{j-1}+(-1)^{a_0}}\otimes v_h)\\
   &&  -Y_{a,h} (i^{-p} u_{a+(-1)^{a_{j-1}}2^{j-1}+(-1)^{a_0}}\otimes v_h)).\\
\end{eqnarray*}

\subsubsection*{Example: Vector fields on $S^{31}$}

In this subsection, we compute explicitly a maximal set of orthogonal linearly independent vector fields on $S^{31}$. Recall that $Cl_{10}^0$ is the biggest even Clifford algebra with $\mathbb{R}^{32}$ as representation space, so there are $9$ linearly independent orthogonal vector fields on the sphere $S^{31}$. Let
$Z\in S^{31}$, in terms of our basis
\begin{eqnarray*}Z&:=&(X_{0} + i Y_{0})u_{0}+ (X_{3} + i Y_{3})u_{3}+ (X_{5} + i Y_{5})u_{5}+
  (X_{6} + i Y_{6})u_{6}+ (X_{9} + i Y_{9})u_{9}+ (X_{10} + i Y_{10})u_{10}\\
  &&+
  (X_{12} + i Y_{12})u_{12}+ (X_{15} + i Y_{15})u_{15}+
  (X_{17} + i Y_{17})u_{17}+ (X_{18} + i Y_{18})u_{18}+
  (X_{20} + i Y_{20})u_{20}\\
  &&
  + (X_{23} + i Y_{23})u_{23}+
  (X_{24} + i Y_{24})u_{24}+ (X_{27} + i Y_{27})u_{27}+
  (X_{29} + i Y_{29})u_{29}+ (X_{30} + i Y_{30})u_{30}
\end{eqnarray*}
then a set of  linearly independent orthogonal vector fields is given by
\begin{eqnarray*}
V_1&=&e_1e_2\cdot Z \\
  &=&(i X_{0} - Y_{0})u_{0}+ (-i X_{3} + Y_{3})u_{3}+
  (-i X_{5} + Y_{5})u_{5}+ (i X_{6} - Y_{6})u_{6}+ (-i X_{9} + Y_{9})u_{9}+
  (i X_{10} - Y_{10})u_{10}\\
  &&+ (i X_{12} - Y_{12})u_{12}+
  (-i X_{15} + Y_{15})u_{15}+ (-i X_{17} + Y_{17})u_{17}+
  (i X_{18} - Y_{18})u_{18}+ (i X_{20} - Y_{20})u_{20}\\
  &&+
  (-i X_{23} + Y_{23})u_{23}+ (i X_{24} - Y_{24})u_{24}+
  (-i X_{27} + Y_{27})u_{27}+ (-i X_{29} + Y_{29})u_{29}+
  (i X_{30} - Y_{30})u_{30}\\
V_2&=&e_1e_3\cdot Z \\
   &=&(-X_{3} - i Y_{3})u_{0}+ (X_{0} + i Y_{0})u_{3}+
  (X_{6} + i Y_{6})u_{5}+ (-X_{5} - i Y_{5})u_{6}+ (X_{10} + i Y_{10})u_{9}+
  (-X_{9} - i Y_{9})u_{10}\\
  &&+ (-X_{15} - i Y_{15})u_{12}+
  (X_{12} + i Y_{12})u_{15}+ (X_{18} + i Y_{18})u_{17}+
  (-X_{17} - i Y_{17})u_{18}+ (-X_{23} - i Y_{23})u_{20}\\
  &&+
  (X_{20} + i Y_{20})u_{23}+ (-X_{27} - i Y_{27})u_{24}+
  (X_{24} + i Y_{24})u_{27}+ (X_{30} + i Y_{30})u_{29}+
  (-X_{29} - i Y_{29})u_{30}\\
V_3&=&e_1e_4\cdot Z\\ 
  &=&(-i X_{3} + Y_{3})u_{0}+ (-i X_{0} + Y_{0})u_{3}+
  (i X_{6} - Y_{6})u_{5}+ (i X_{5} - Y_{5})u_{6}+ (i X_{10} - Y_{10})u_{9}+
  (i X_{9} - Y_{9})u_{10}\\
  &&+ (-i X_{15} + Y_{15})u_{12}+
  (-i X_{12} + Y_{12})u_{15}+ (i X_{18} - Y_{18})u_{17}+
  (i X_{17} - Y_{17})u_{18}+ (-i X_{23} + Y_{23})u_{20}\\
  &&+
  (-i X_{20} + Y_{20})u_{23}+ (-i X_{27} + Y_{27})u_{24}+
  (-i X_{24} + Y_{24})u_{27}+ (i X_{30} - Y_{30})u_{29}+
  (i X_{29} - Y_{29})u_{30}\\
V_4&=&e_1e_5\cdot Z\\ 
   &=&(X_{5} + i Y_{5})u_{0}+ (X_{6} + i Y_{6})u_{3}+
         (-X_{0} - i Y_{0})u_{5}+ (-X_{3} - i Y_{3})u_{6}+
         (-X_{12} - i Y_{12})u_{9}+ (-X_{15} - i Y_{15})u_{10}\\
         &&+
         (X_{9} + i Y_{9})u_{12}+ (X_{10} + i Y_{10})u_{15}+
         (-X_{20} - i Y_{20})u_{17}+ (-X_{23} - i Y_{23})u_{18}+
         (X_{17} + i Y_{17})u_{20}\\
         &&+ (X_{18} + i Y_{18})u_{23}+
         (X_{29} + i Y_{29})u_{24}+ (X_{30} + i Y_{30})u_{27}+
         (-X_{24} - i Y_{24})u_{29}+ (-X_{27} - i Y_{27})u_{30}\\
V_5&=&e_1e_6\cdot Z \\
  &=&(i X_{5} - Y_{5})u_{0}+ (i X_{6} - Y_{6})u_{3}+
  (i X_{0} - Y_{0})u_{5}+ (i X_{3} - Y_{3})u_{6}+ (-i X_{12} + Y_{12})u_{9}+
  (-i X_{15} + Y_{15})u_{10}\\
  &&+ (-i X_{9} + Y_{9})u_{12}+
  (-i X_{10} + Y_{10})u_{15}+ (-i X_{20} + Y_{20})u_{17}+
  (-i X_{23} + Y_{23})u_{18}+ (-i X_{17} + Y_{17})u_{20}\\
  &&+
  (-i X_{18} + Y_{18})u_{23}+ (i X_{29} - Y_{29})u_{24}+
  (i X_{30} - Y_{30})u_{27}+ (i X_{24} - Y_{24})u_{29}+
  (i X_{27} - Y_{27})u_{30}\\
V_6&=&e_1e_7\cdot Z \\
  &=&(-X_{9} - i Y_{9})u_{0}+ (-X_{10} - i Y_{10})u_{3}+
  (-X_{12} - i Y_{12})u_{5}+ (-X_{15} - i Y_{15})u_{6}+
  (X_{0} + i Y_{0})u_{9}+ (X_{3} + i Y_{3})u_{10}\\
  &&+ (X_{5} + i Y_{5})u_{12}+
  (X_{6} + i Y_{6})u_{15}+ (X_{24} + i Y_{24})u_{17}+
  (X_{27} + i Y_{27})u_{18}+ (X_{29} + i Y_{29})u_{20}\\
  &&+
  (X_{30} + i Y_{30})u_{23}+ (-X_{17} - i Y_{17})u_{24}+
  (-X_{18} - i Y_{18})u_{27}+ (-X_{20} - i Y_{20})u_{29}+
  (-X_{23} - i Y_{23})u_{30}\\
V_7&=&e_1e_8\cdot Z \\
   &=&(-i X_{9} + Y_{9})u_{0}+ (-i X_{10} + Y_{10})u_{3}+
       (-i X_{12} + Y_{12})u_{5}+ (-i X_{15} + Y_{15})u_{6}+
       (-i X_{0} + Y_{0})u_{9}+ (-i X_{3} + Y_{3})u_{10}\\
       &&+
       (-i X_{5} + Y_{5})u_{12}+ (-i X_{6} + Y_{6})u_{15}+
       (i X_{24} - Y_{24})u_{17}+ (i X_{27} - Y_{27})u_{18}+
       (i X_{29} - Y_{29})u_{20}\\
       &&+ (i X_{30} - Y_{30})u_{23}+
       (i X_{17} - Y_{17})u_{24}+ (i X_{18} - Y_{18})u_{27}+
       (i X_{20} - Y_{20})u_{29}+ (i X_{23} - Y_{23})u_{30}\\
V_8&=&e_1e_9\cdot Z \\
   &=&(X_{17} + i Y_{17})u_{0}+ (X_{18} + i Y_{18})u_{3}+
      (X_{20} + i Y_{20})u_{5}+ (X_{23} + i Y_{23})u_{6}+
      (X_{24} + i Y_{24})u_{9}+ (X_{27} + i Y_{27})u_{10}\\
      &&+
      (X_{29} + i Y_{29})u_{12}+ (X_{30} + i Y_{30})u_{15}+
      (-X_{0} - i Y_{0})u_{17}+ (-X_{3} - i Y_{3})u_{18}+
      (-X_{5} - i Y_{5})u_{20}\\
      &&+ (-X_{6} - i Y_{6})u_{23}+
      (-X_{9} - i Y_{9})u_{24}+ (-X_{10} - i Y_{10})u_{27}+
      (-X_{12} - i Y_{12})u_{29}+ (-X_{15} - i Y_{15})u_{30}\\
V_9&=&e_1e_{10}\cdot Z \\
  &=&(i X_{17} - Y_{17})u_{0}+ (i X_{18} - Y_{18})u_{3}+
 (i X_{20} - Y_{20})u_{5}+ (i X_{23} - Y_{23})u_{6}+
 (i X_{24} - Y_{24})u_{9}+ (i X_{27} - Y_{27})u_{10}\\
 &&+
 (i X_{29} - Y_{29})u_{12}+ (i X_{30} - Y_{30})u_{15}+
 (i X_{0} - Y_{0})u_{17}+ (i X_{3} - Y_{3})u_{18}+ (i X_{5} - Y_{5})u_{20}\\
 &&+
 (i X_{6} - Y_{6})u_{23}+ (i X_{9} - Y_{9})u_{24}+
 (i X_{10} - Y_{10})u_{27}+ (i X_{12} - Y_{12})u_{29}+
 (i X_{15} - Y_{15})u_{30}.
\end{eqnarray*}
In terms of coordinate vectors, one can write these vector fields as follows
\begin{eqnarray*}
V_1 &=&
(-Y_{0},X_{0},Y_{3},-X_{3},Y_{5},-X_{5},-Y_{6},X_{6},Y_{9},-X_{9},-Y_{10},X_{10},-Y_{12},X_{12},Y_{15},\\&&-X_{15},Y_{17},-X_{17},-Y_{18},X_{18},-Y_{20},X_{20},Y_{23},-X_{23},-Y_{24},X_{24},Y_{27},-X_{27},Y_{29},-X_{29},-Y_{30},X_{30}),
\\
V_2&=&
(-X_{3},-Y_{3},X_{0},Y_{0},X_{6},Y_{6},-X_{5},-Y_{5},X_{10},Y_{10},-X_{9},-Y_{9},-X_{15},-Y_{15},X_{12},Y_{12},\\&&X_{18},Y_{18},-X_{17},-Y_{17},-X_{23},-Y_{23},X_{20},Y_{20},-X_{27},-Y_{27},X_{24},Y_{24},X_{30},Y_{30},-X_{29},-Y_{29}),
\\
V_3&=&
(Y_{3},-X_{3},Y_{0},-X_{0},-Y_{6},X_{6},-Y_{5},X_{5},-Y_{10},X_{10},-Y_{9},X_{9},Y_{15},-X_{15},Y_{12},-X_{12},\\&&-Y_{18},X_{18},-Y_{17},X_{17},Y_{23},-X_{23},Y_{20},-X_{20},Y_{27},-X_{27},Y_{24},-X_{24},-Y_{30},X_{30},-Y_{29},X_{29}),
\\
V_4 &=&(X_{5},Y_{5},X_{6},Y_{6},-X_{0},-Y_{0},-X_{3},-Y_{3},-X_{12},-Y_{12},-X_{15},-Y_{15},X_{9},Y_{9},X_{10},Y_{10},\\&&-X_{20},-Y_{20},-X_{23},-Y_{23},X_{17},Y_{17},X_{18},Y_{18},X_{29},Y_{29},X_{30},Y_{30},-X_{24},-Y_{24},-X_{27},-Y_{27}),
\\
V_5 &=&(-Y_{5},X_{5},-Y_{6},X_{6},-Y_{0},-X_{0},-Y_{3},X_{3},Y_{12},-X_{12},Y_{15},-X_{15},Y_{9},-X_{9},Y_{10},-X_{10},\\&&Y_{20},-X_{20},Y_{23},-X_{23},Y_{17},-X_{17},Y_{18},-X_{18},-Y_{29},X_{29},-Y_{30},X_{30},-Y_{24},X_{24},-Y_{27},X_{27})),
\\
V_6&=&
(-X_{9},-Y_{9},-X_{10},-Y_{10},-X_{12},-Y_{12},-X_{15}, Y_{15},X_{0},Y_{0},X_{3},Y_{3},X_{5},Y_{5},X_{6},Y_{6},\\&&X_{24},Y_{24},X_{27},Y_{27},X_{29},Y_{29},X_{30},Y_{30},-X_{17},-Y_{17},-X_{18},-Y_{18},-X_{20},-Y_{20},-X_{23},-Y_{23}),
\\
V_7 &=&(Y_{9},-X_{9},Y_{10},-X_{10},Y_{12},-X_{12},Y_{15},-X_{15},Y_{0},-X_{0},Y_{3},-X_{3},Y_{5},-X_{5},Y_{6},-X_{6},\\&&-Y_{24},X_{24},-Y_{27},X_{27},-Y_{29},X_{29},-Y_{30},X_{30},-Y_{17},X_{17},-Y_{18},X_{18},-Y_{20},X_{20},-Y_{23},X_{23}),
\\
V_8 &=&(X_{17},Y_{17},X_{18},Y_{18},X_{20},Y_{20},X_{23},Y_{23},X_{24},Y_{24},X_{27},Y_{27},X_{29},Y_{29},X_{30},Y_{30},\\&& -X_{0},-Y_{0},-X_{3},-Y_{3},-X_{5},-Y_{5},-X_{6},-Y_{6},-X_{9},-Y_{9},-X_{10},-Y_{10},-X_{12},-Y_{12},-X_{15},-Y_{15}),
\\
V_9 &=&
(-Y_{17},X_{17},-Y_{18},X_{18},-Y_{20},X_{20},-Y_{23},X_{23},-Y_{24},X_{24},-Y_{27},X_{27},-Y_{29},X_{29},-Y_{30},X_{30},\\&&-Y_{0}, X_{0},-Y_{3},X_{3},-Y_{5},X_{5},-Y_{6},X_{6},-Y_{9},X_{9},-Y_{10},X_{10},-Y_{12},X_{12},-Y_{15},X_{15}).
\end{eqnarray*}
Relabelling the entries of $Z$, 
\begin{eqnarray*}Z&:=&(v_1,v_2,v_3,v_4,v_5,v_6,v_7,v_8,v_9,v_{10},v_{11},v_{12},v_{13},v_{14},v_{15},v_{16},\\&& v_{17},v_{18},v_{19},
v_{20},v_{21},v_{22},v_{23},v_{24},v_{25},v_{26},v_{27},v_{28},v_{29},v_{30},v_{31},v_{32}),
\end{eqnarray*}
the vector fields are the following
\begin{eqnarray*}
V_1&=&
(-v_2,v_1,v_4,-v_3,v_6,-v_5,-v_8,v_7,v_{10},-v_9,-v_{12},v_{11},-v_{14},v_{13},v_{16},-v_{15},\\&&v_{18},-v_{17},-v_{20},v_{19},-v_{22},v_{21},v_{24},-v_{23},-v_{26},v_{25},v_{28},-v_{27},v_{30},-v_{29},-v_{32},v_{31}),
\\
V_2&=&
(-v_3,-v_4,v_1,v_2,v_7,v_8,-v_5,-v_6,v_{11},v_{12},-v_9,-v_{10},-v_{15},-v_{16},v_{13},v_{14},\\&&v_{19},v_{20},-v_{17},-v_{18},-v_{23},-v_{24},v_{21},v_{22},-v_{27},-v_{28},v_{25},v_{26},v_{31},v_{32},-v_{29},-v_{30}),
\\
V_3 &=&
(v_4,-v_3,v_2,-v_1,-v_8,v_7,-v_6,v_5,-v_{12},v_{11},-v_{10},v_9,v_{16},-v_{15},v_{14},-v_{13},\\&&-v_{20},v_{19},-v_{18},v_{17},v_{24},-v_{23},v_{22},-v_{21},v_{28},-v_{27},v_{26},-v_{25},-v_{32},v_{31},-v_{30},v_{29}),
\\
V_4 &=&(v_5,v_6,v_7,v_8,-v_1,-v_2,-v_3,-v_4,-v_{13},-v_{14},-v_{15},-v_{16},v_9,v_{10},v_{11},v_{12},\\&&-v_{21},-v_{22},-v_{23},-v_{24},v_{17},v_{18},v_{19},v_{20},v_{29},v_{30},v_{31},v_{32},-v_{25},-v_{26},-v_{27},-v_{28}),
\\
V_5 &=&(-v_6,v_5,-v_8,v_7,-v_2,-v_1,-v_4,v_3,v_{14},-v_{13},v_{16},-v_{15},v_{10},-v_9,v_{12},-v_{11},\\&&v_{22},-v_{21},v_{24},-v_{23},v_{18},-v_{17},v_{20},-v_{19},-v_{30},v_{29},-v_{32},v_{31},-v_{26},v_{25},-v_{28},v_{27}),
\\
V_6 &=&
(-v_9,-v_{10},-v_{11},-v_{12},-v_{13},-v_{14},-v_{15}, v_{16},v_1,v_2,v_3,v_4,v_5,v_6,v_7,v_8,\\
&&v_{25},v_{26},v_{27},v_{28},v_{29},v_{30},v_{31},v_{32},-v_{17},-v_{18},-v_{19},-v_{20},-v_{21},-v_{22},-v_{23},-v_{24}),
\\
V_7 &=&(v_{10},-v_9,v_{12},-v_{11},v_{14},-v_{13},v_{16},-v_{15},v_2,-v_1,v_4,-v_3,v_6,-v_5,v_8,-v_7,\\&&-v_{26},v_{25},-v_{28},v_{27},-v_{30},v_{29},-v_{32},v_{31},-v_{18},v_{17},-v_{20},v_{19},-v_{22},v_{21},-v_{24},v_{23}),
\\
V_8 &=&(v_{17},v_{18},v_{19},v_{20},v_{21},v_{22},v_{23},v_{24},v_{25},v_{26},v_{27},v_{28},v_{29},v_{30},v_{31},v_{32},\\&& -v_1,-v_2,-v_3,-v_4,-v_5,-v_6,-v_7,-v_8,-v_9,-v_{10},-v_{11},-v_{12},-v_{13},-v_{14},-v_{15},-v_{16}),
\\
V_9&=&
(-v_{18},v_{17},-v_{20},v_{19},-v_{22},v_{21},-v_{24},v_{23},-v_{26},v_{25},-v_{28},v_{27},-v_{30},v_{29},-v_{32},v_{31},\\&&-v_2, v_1,-v_4,v_3,-v_6,v_5,-v_8,v_7,-v_{10},v_9,-v_{12},v_{11},-v_{14},v_{13},-v_{16},v_{15}).
\end{eqnarray*}

{\small 
\renewcommand{\baselinestretch}{0.5}
\newcommand{\bi}{\vspace{-.05in}\bibitem} }

\enddocument